\documentclass{article}
\usepackage{amsfonts,amsmath,theorem}
\newtheorem{lemma}{Lemma}[section]
\newtheorem{theorem}{Theorem}[section]

\newtheorem{remark}{Remark}[section]
\newtheorem{prop}{Proposition}[section]
\newtheorem{proposition}{Proposition}[section]
\newtheorem{defi}{Definition}[section]

\newcommand{\nn}{\nonumber}

\newcommand{\be}{\begin{equation}}
\newcommand{\ee}{\end{equation}}
\newcommand{\bea}{\begin{eqnarray}}
\newcommand{\eea}{\end{eqnarray}}
\newcommand{\bd}{\begin{displaymath}}
\newcommand{\ed}{\end{displaymath}}

\newcommand{\U}{{\mathcal U}}

\def \R{{\mathbf{R}}}

\def \UU{{{\mathcal U}}}

\def\ZZ {{{\mathbf Z}}}

\def \Re{Re}                    

\def \eps{{\epsilon}}

\def \uT {\tilde u}

\def \del {{\partial}}

\title{Euler Equations of Incompressible Ideal Fluids}

\author{
Claude BARDOS%
\footnote{ Universit\'e Denis Diderot  and Laboratory JLL
Universit\'e Pierre et Marie Curie, Paris, France
(bardos@ann.jussieu.fr).}, and
Edriss S. TITI%
\footnote{ Department of Mathematics and  Department of Mechanical
and  Aerospace Engineering University of California Irvine, CA
92697-3875, USA. {\bf Also}, Department of Computer Science and
Applied Mathematics, Weizmann Institute of Science Rehovot 76100,
Israel
(etiti@math.uci.edu).}}%

\date{February 27, 2007 }

\begin{document}
\maketitle
\begin{abstract}
This article is
a survey concerning the state-of-the-art mathematical theory of the
Euler equations of incompressible homogenous ideal fluid. Emphasis
is put on the different types of emerging instability, and how they
may be related to the description of turbulence.
\end{abstract}
\section{Introduction}


This contribution is mostly devoted to the time dependent analysis
of the $2d$ and $3d$  Euler equations
\begin{equation}
\del_t u + \nabla\cdot (u \otimes u) +\nabla p=0\,,\,\, \nabla\cdot
u =0\,,\label{Euler0}
\end{equation}
of incompressible homogenous ideal fluid. We intend to connect
several known (and maybe less known) points of view concerning this
very classical problem. Furthermore, we will investigate the
conditions under which one can consider the above problem as the
limit of the incompressible Navier--Stokes equations:
\begin{equation}
\del_t u_\nu + \nabla \cdot  (u_\nu  \otimes u_\nu)-\nu\Delta u_\nu
+\nabla p_\nu=0\,,\,\,\, \nabla\cdot u_\nu
=0\,,\label{NS1}
\end{equation}
when the viscosity $\nu \to 0$, i.e.,~as the Reynolds number goes to
infinity.

 At the macroscopic level the  Reynolds number, $\Re$,
  corresponds to the ratio of {\it the strength of the
  nonlinear effects} and {\it the strength of the linear viscous
  effects}. Therefore, with the introduction
 of  a characteristic velocity, $U$, and  a characteristic
 length scale, $L$,  of the flow one has the dimensionless
 parameter:
\begin{equation} \Re=\frac{UL}{\nu} \label{Re}\,.
\end{equation}
With the introduction of the characteristic time scale $T=L/U$ and
the dimensionless variables:
\begin{equation*}
x'=\frac{x}L\,, t'=\frac{t}{T} \,\hbox{ and  } u'=\frac {u'} {U}\,,
\end{equation*}
the Navier--Stokes equations (\ref{NS1}) take  the non-dimensional
form:
 \begin{equation}
\del_t u' + \nabla_{x'} \cdot (u' \otimes u')-\frac 1{\Re} \Delta_{x'} u' +\nabla_{x'} p'=0\,,\,\, \nabla\cdot u' =0\,.\label{NSAD}
\end{equation}
These are the equations to be considered in the sequel, omitting the
$'$ and returning to the notation $\nu$ for $\Re^{-1}\,.$

In the presence of  physical boundary the problems (\ref{Euler0})
and (\ref{NS1}) will be considered in the open domain $\Omega\subset
\R^d\,, d=2\,,\,\, d=3\,,$ with a piecewise smooth boundary $\del
\Omega\,.$

There are several good reason to focus at present on the
``mathematical analysis" of the Euler equations rather than on the
Navier--Stokes equations.

1. Turbulence applications involving the Navier--Stokes equations
(\ref{NSAD}) often correspond to very large Reynolds numbers; and a
theorem which is valid for any finite, but very large, Reynolds
number is expected to be compatible with results concerning infinite
Reynolds numbers. In fact, this is the case when $\Re=\infty$ which
drives other results and we  will give several examples of this
fact.

2.  Many nontrivial and sharp results obtained for the
incompressible Navier--Stokes equations rely  on the smoothing
effect of the Laplacian,  with viscosity $\nu>0$, and on the
invariance of the set of solutions under the scaling:
\begin{equation}
u(x,t) \mapsto \lambda u(\lambda x,\lambda^2 t)\,.\label{scaling}
\end{equation}
However, simple examples with the same scalings, but without  a
conservation law of energy may exhibit very different behavior
concerning regularity and stability.

1. With $\phi$ being a scalar function, the viscous Hamilton--Jacobi
type or Burgers equation
\begin{eqnarray}
&&\del_t \phi -\nu\Delta \phi+\frac12|\nabla \phi|^2=0\,
\hbox{ in } \Omega\times \R^+_t\,, \label{Bur1}\\
&&\phi(x,t)=0 \hbox{ for } x\in\del \Omega\,, \hbox{ and }
\phi(\cdot,0)=\phi_0(\cdot)\in L^\infty(\Omega)\,, \nonumber
\end{eqnarray}
has (because of the maximum principle) a global smooth solution, for
$\nu>0$. However, for $\nu=0$, it is well known that certain
solutions of the inviscid Burgers  equation (\ref{Bur1}) will become
singular (with shocks) in finite time.

2. Denote by $|\nabla|$ the square root of the operator $-\Delta$,
defined in $\Omega$ with  Dirichlet homogeneous boundary conditions.
Consider the solution $u(x,t)$ of the equation
\begin{eqnarray}
&&\del_t u -\nu\Delta u+\frac12|\nabla |(u^2)=0\,\, \hbox{ in }\,  \label{cheap1}\Omega\times \R^+_t\,,\\
&&u(x,t)=0 \hbox{ for } x\in\del \Omega\,, \hbox{ and }
u(\cdot,0)=u_0(\cdot)\in L^\infty(\Omega)\,.\label{cheap2}
\end{eqnarray}
Then one has the following proposition.
\begin{prop}
Assume that the initial data $u_0$ satisfies the relation:
\begin{equation}
\int_\Omega u_0(x)\phi_1(x)dx =-M <0\,,
\end{equation}
where $\phi_1(x)\ge 0$ denotes the first eigenfunction of the
operator $-\Delta$ (with Dirichlet boundary condition),
$-\Delta\phi_1=\lambda_1\phi_1\,.$  Then if $M$ is large enough, the
corresponding solution $u(x,t)$ of the system (\ref{cheap1})
(\ref{cheap2}) blows up in a finite time.
\end{prop}
{\bf Proof.} The $L^2$ scalar product of the equation (\ref{cheap1})
with $\phi_1(x)$ gives
\begin{eqnarray*}
&&\frac{d}{dt}\int_\Omega u(x,t)\phi_1(x)dx +\nu \lambda_1\int_\Omega u(x,t)\phi_1(x)dx\\
&&\hskip 1in =-\frac{\sqrt{\lambda_1} }2\int_\Omega
u(x,t)^2\phi_1(x)dx\,.
\end{eqnarray*}
Since $\phi_1(x) \ge 0$  then the Cauchy--Schwarz inequality implies
$$
\left(\int u(x,t)\phi_1(x)dx\right)^2\le \left(\int_\Omega
u(x,t)^2\phi_1(x)dx\right)\left(\int_\Omega \phi_1(x)dx\right)\,.
$$
As a result of the above the quantity $m(t)=-\int_\Omega
u(x,t)\phi_1(x)dx $ satisfies the relation:
\begin{equation*}
\frac{dm}{dt} +\lambda_1m\ge \frac{ {\sqrt{\lambda_1}\int_\Omega
\phi_1(x)dx}  } 2m^2\,,\, {\hbox {with}}\, m(0)=M\,;
\end{equation*}
and the conclusion of the proposition follows.
\begin{remark}
The above example has been introduced with $\Omega=\R^3$ by
Montgomery--Smith \cite{MS} under the name ``cheap Navier--Stokes
equations" with the purpose of underlying the role of the
conservation of energy (which is not present in the above examples)
in the Navier--Stokes dynamics. His proof shows that the same blow
up property may appear in any space dimension for the solution of
the ``cheap hyper-viscosity  equations"
\begin{equation*}
\del_t u +\nu(-\Delta)^m u+\frac12|\nabla | (u^2)=0\,.
\end{equation*}
On the other hand, one should observe that the  above argument does
not apply to the Kuramoto--Sivashinsky-like equations
\begin{equation}
\del_t \phi +\nu(-\Delta)^m \phi + \alpha\Delta \phi +\frac12|\nabla
\phi|^2=0 \label{bbt},
\end{equation}
for $m\ge 2$.  Without a maximum principle or without the control of
some sort of  energy the question of global existence of smooth
solution,  or finite time blow up of some solution, to the above
equation is an open problem in $\R^n$, for $ n \ge 2$ and for $m\ge
2$. However, if in (\ref{bbt}) the term $|\nabla \phi|^2$ is
replaced by $ |\nabla \phi|^{2+\gamma}\,, \gamma>0$ one can prove
the blow up of some solutions (cf. \cite{BBT} and references
therein).
\end{remark}

In conclusion, the above examples indicate  that the conservation of
some sort of  energy, which is guaranteed by the structure of the
equation, is essential in the analysis of the dynamics of the
underlying problem. In particular, this very basic fact plays an
essential role  in the dynamics of the Euler equations.

Taking into account the above simple examples, the rest of the paper
is organized as follows. In  section 2 classical existence and
regularity results for the time dependent Euler equations are
presented. Section 3 provides more examples concerning the
pathological behavior of solutions of the Euler equations. The fact
that the solutions of the Euler equations may exhibit oscillatory
behavior  implies similar behavior for the solutions of the
Navier--Stokes equations, as the viscosity tends to zero. The
existence of (or lack thereof) strong convergence is analyzed in
section 4 with the introduction of the Reynolds stresses tensor, and
the notion of {\it dissipative solution}.  A standard and very
important problem, for both theoretical study  and applications, is
the vanishing viscosity limit of solutions of the Navier--Stokes
equations subject to the no-slip Dirichlet boundary condition, in
domains with physical boundaries. Very few mathematical results are
available for this very unstable situation. One of the most striking
results is a theorem of Kato \cite{KA}, which is presented in
section 5. Section 6 is again devoted to the Reynolds stresses
tensor. We show that with the introduction of the Wigner measure the
notion of Reynolds stresses tensor, deduced from the defect in
strong convergence as the viscosity tends to zero, plays the same
role as the one originally introduced in the statistical theory of
turbulence. When the zero viscosity limit of solutions of the
Navier--Stokes equations  agrees with the solution of the Euler
equations the main difference is confined in a boundary layer which
is described by the Prandtl equations. These equations are briefly
described in section 7. It is also recalled how the mathematical
results are in agreement with the instability of the physical
problem. The Kelvin--Helmholtz problem exhibits also some basic
similar instabilities, but it is in some sense simpler. This is
explained at the end of section 7, where it is also shown that some
recent results of \cite{LE} ,  \cite{Wu1} and \cite{Wu2}, on the
regularity of the vortex sheet (interface), do contribute to the
understanding of the instabilities of the original problem.

\section{Classical existence and regularity results }

\subsection{Introduction}

The Euler equations correspond, formally, to the limit case when the
viscosity is $0$, or the Reynolds number is infinite:
\begin{equation}
\del_t u + \nabla \cdot (u \otimes u) + \nabla p=0\, ,\,\,
\nabla\cdot u =0\,, \hbox{ in } \Omega \,.\label{Euler1}
\end{equation}
In the presence of physical boundaries, the above system is
supplemented with the standard, no-normal flow, boundary condition:
\begin{equation}
u\cdot \vec{n}=0\,\, \hbox{ on } \del\Omega\,, \label{Euler-BC}
\end{equation}
where $\vec{n}$ denotes the  outwards normal vector to the boundary
$\del\Omega$. It turns out that the vorticity, $\omega=\nabla \wedge
u$, is ``the basic quantity", from both  the physical and
mathematical analysis points of view. Therefore,  equations
(\ref{Euler1}) and (\ref{Euler-BC}), written in terms of the
vorticity, are equivalent to the system:
\begin{eqnarray}
&&\del_t \omega + u\cdot \nabla \omega =\omega \cdot \nabla u\,\, {\hbox{in}}\, \Omega \label{vort}\\
&&\nabla \cdot u=0\,, \nabla \wedge u=\omega\,\, {\hbox{in}}\,
\Omega\,,\,\, {\hbox{and}}\, u\cdot \vec n =0\,\hbox{ on }
\del\Omega\,. \label{ell2}
\end{eqnarray}
 That is,  system (\ref{ell2}) fully determines $u$ in terms
 of $\omega$, which makes the above system ``closed".
 More precisely, the operator $K:\omega \mapsto u$
 defined by  relation (\ref{ell2}) is a linear continuous
 map from $C^\alpha(\Omega)$ with values in
 $C^{\alpha +1}(\Omega)$ (with $\alpha>0$); and from
 $H^s(\Omega)$ with values in $H^{s+1}(\Omega)\,.$

 Furthermore, for $2d$ flows, the vorticity is perpendicular
 to the plane of motion and therefore  equation
 (\ref{vort}) is reduced (this can also be checked directly)
 to the advection equation
 \begin{equation}
 \del_t \omega +u\cdot \nabla \omega =0\,.
 \end{equation}
 The structure of the quadratic nonlinearity in  (\ref{vort}) has the
 following consequences, which are described below. We will be
 presenting only the essence of the
 essential arguments and not the full details of the proofs (see,
e.g., \cite{Majda-Bertozzi} or \cite{MP} for the details).
 \subsection{General results in $3d$}

 The short time existence of a smooth solution for
 the $3d$ incompressible Euler equations has been
 obtained already a long time ago, provided the initial
 data are smooth enough.  To the best of our knowledge the
 original proof goes back to Lichtenstein \cite{LI}.
 The proof is based on a nonlinear Gronwall estimate
 of the following type:
\begin{equation}
y'\le  C y^{\frac 32} \Rightarrow y(t)\le \frac{y(0)}{(1-2t C
y^{\frac12}(0))^2}\,. \label{Gronwall}
\end{equation}
Therefore, the value of $y(t)$, which represents an adequate norm of
the solution, is finite for a finite interval of time; which depends
on the size of the initial value of $y(0)$, i.e. the initial data of
the solution of Euler. These initial data have to be chosen from an
appropriate space of regular enough functions. In particular, if we
consider the solution in the Sobolev space $H^s$, with
$s>\frac{5}{2}$, then by taking the scalar product, in the Sobolev
space $H^s$, of the Euler equations with the solution $u$,  and  by
using the appropriate Sobolev estimates we obtain:
\begin{equation}
\frac12\frac{d |\!|u |\!|_{H^s}^2}{dt}= -(\nabla\cdot (u\otimes u),
u)_{H^s} \le C_s |\!|u |\!|_{H^s}^2|\!| \nabla u|\!|_{L^\infty} \le
C|\!|u |\!|_{H^s}^3\, . \label{scalarprod}
\end{equation}
As a result of (\ref{Gronwall}) and (\ref{scalarprod})  we obtain
the local, in time, existence of a smooth solution.

As  standard in many nonlinear time dependent problems local
regularity of smooth (strong) solutions implies local uniqueness and
local stability (i.e, continuous dependence on initial data).
Furthermore, one may exhibit a threshold for this existence,
uniqueness, and propagation of the regularity of the intial data
(including analyticity Bardos and Benachour \cite{BB}). More
precisely one uses the following theorem.
\begin{theorem}{\bf Beale--Kato--Majda} \cite{BKM}
Let $u(t)$ be a solution of the $3d$ incompressible Euler equations
which  is regular for $0\le t  <T$; that is,
\begin{equation*}
{\hbox{for all}}\,\,\, t \in [0, T]\,,\,\,\,   u(t)\in
H^{s}(\Omega)\,,\, {\hbox{for some}}\,\,\, s>\frac53\,.
\end{equation*}
Assume that
\begin{equation}
\int_0^T|\!|\nabla \wedge u(.,t)|\!|_{L^\infty} dt <\infty\,,
\label{bmk}
\end{equation}
then $u(t)$  can be uniquely extended up to a time $T+\delta$
($\delta>0$) as a smooth solution of the Euler equations.
\end{theorem}
The main interest of this statement is the fact  that it shows that
if one starts with smooth initial data, then  instabilities appears
only if the size of the   vorticity becomes arbitrary  large.
\begin{remark}
The Beale--Kato--Majda theorem was first proven in the whole space
in \cite{BKM}. Extension to a periodic ``box" is easy. For a bounded
domain with the boundary condition $u\cdot \vec n =0$ it was
established by Ferrari \cite{FE}. By combining arguments form
\cite{BB} and \cite{FE} one can show, as in the  Beale--Kato--Majda
theorem, that the solution of $3d$ Euler equations, with real
analytic initial data, remains real analytic  as long as (\ref{bmk})
holds.
\end{remark}
The Beale--Kato--Majda result has been slightly improved by Kozono
\cite{KO} who proved that on the left-hand side of (\ref{bmk}), the
$|\!|\cdot|\!|_{L^\infty}$--norm can be replaced by the norm in the
BMO space. This generalization is interesting because it adapts
harmonic analysis (or Fourier modes decomposition) techniques which
is an important tool for the study of ``turbulent" solutions;
indeed, the space BMO, as the dual space of the Hardy space
${\mathcal H}^1$\,, is well defined in the frequency (Fourier)
space. In fact, cf. \cite{ME}, BMO is the smallest space containing
$L^\infty$\,, which is also invariant under the action of a zero
order pseudodifferential operators. The idea behind the
Beale--Kato--Majda theorem, and its improvement, is the fact the
solution $u$ of  the elliptic equations (\ref{ell2}) satisfies the
following estimate, for $1<p<\infty$\,,
\begin{equation}
|\!| \nabla u |\!|_{W^{s,p}}\le C_{s,p}\left(|\!| u |\!|_{W^{s ,p}}
+ |\!| \omega |\!|_{W^{s ,p}}\right)\,. \label{el3}
\end{equation}
This relation could also be phrased in the context of H\"older
spaces $C^{k,\alpha}\,,\alpha>0$. We stress, however, that the
estimate (\ref{el3}) ceases to be true for $p=\infty$ (or
$\alpha=0\,.$) This is  due to the nature of the singularity,  which
is of the form $|x-y|^{2-d}$\,,  in the kernel of the operator $K$
and which leads (for  $s>d/2+1$) to the estimate:
\begin{eqnarray}
&&|\!| \nabla u |\!|_{L^\infty }
\le C\left(  |\!| \omega |\!|_{L^\infty } \log(1 + |\!| u |\!|^2_{ H^s} )\right)\,, \label{loglip}\\
&&\hbox { or sharper } \,|\!| \nabla u |\!|_{L^\infty }\le C\left(
|\!| \omega |\!|_{BMO } \log(1 + |\!| u |\!|^2_{ H^s} )\right)
\,.\label{loglip2}
\end{eqnarray}
With $z=1 + |\!| u |\!|^2_{ H^s} $ and thanks to (\ref{loglip2}) the
inequality (\ref{scalarprod}) becomes:
\begin{equation*}
\frac{d}{dt}z\le C |\!| \omega |\!|_{BMO}\, z\log z \,.
\end{equation*}
This yields:
\begin{equation*}
(1 + |\!| u(t) |\!|^2_{ H^s}  )\le (1 + |\!|u(0)|\!|^2_{ H^s} )^{
e^{C\int_0^t|\!| \omega (s)|\!|_{BMO }ds} }\,,
\end{equation*}
which proves the statement. The uniqueness of solutions can be
proven along the same lines, as long as
$$
\int_0^t |\!| \omega (s) |\!|_{BMO } $$
remains finite.
\begin{remark}
The vorticity $\omega$ can be represented by  the anti-symmetric
part of the deformation tensor $\nabla u\,. $ However, in the
estimates (\ref{loglip}) or (\ref{loglip2})  this anti-symmetric
part, i.e. $\omega$, can be replaced by the symmetric part of the
deformation tensor
\begin{equation}S(u)=\frac12(\nabla u +(\nabla u)^t)\,.
\end{equation}
Therefore, the theorems of Beale--Kato--Majda and Kozono can be
rephrased in terms of this symmetric tensor \cite{KO}.
\end{remark}
In fact, the above deformation tensor  $S(u)$ (or $\tilde S(\omega)$
when expressed in term of the vorticity), plays an important role in
a complementary result of Constantin, Fefferman and Majda
\cite{CFM}, which shows that it is mostly the variations in the
direction of the vorticity that may produce singularities.
\begin{prop}
\label{ Constantin Fefferman Majda} \cite{CFM} Let $u$, which is
defined  in $Q= \Omega\times(0,T)$, be a smooth solution of the
Euler equations. Introduce the quantities $k_1(t)$ and $k_2(t)$
(which are well defined for $t<T$):
\begin{equation*}
k_1(t)=\sup_{x\in \Omega} |u(x,t)|\,,
\end{equation*}
which measures the size of the velocity, and
\begin{equation*}
k_2(t)= 4 \pi \sup_{x,y\in\Omega\,,\, x\not=y}
\frac{|\xi(x,t)-\xi(y,t)|}{|x-y|}\,,
\end{equation*}
which measures the Lipschitz regularity of the direction
$$
\xi(x,t)=\frac{\omega(x,t)}{|\omega(x,t)|}
$$
of the vorticity. Then under the assumptions
\begin{equation}
\int_0^T(k_1(t)+k_2(t))dt<\infty\,\,\, {\hbox{and}}\,\, \int_0^T
k_1(t)k_2(t)dt<\infty \,,
\end{equation}
the solution $u$ exists, and is as smooth as the initial data up to a
time $T+\delta$ for some $\delta>0\,.$
\end{prop}

{\bf Proof.} As before we only present here the basic ideas, and for simplicity we will focus on the case when
$\Omega=\R^3$. First, since
\begin{equation}
S(u)=\frac12\left(\nabla u +(\nabla u)^t\right)(x,t)=\tilde
S(\omega)(x,t)\,,
\end{equation}
we have
\begin{equation}
\frac12(\del_t|\omega|^2+ u\cdot \nabla |\omega|^2)=(\omega\cdot \nabla u, \omega)=(\tilde S(\omega)\omega , \omega)\,,
\end{equation}
which gives
\begin{equation}
\frac{d|\!|\omega|\!|_{\infty}}{dt}\le \sup_x(|\tilde
S(\omega)|)|\!|\omega|\!|_{\infty}\,.\label{gron}
 \end{equation}
Next, we consider only the singular part of the operator
$\omega\mapsto\tilde S(\omega)\,.$ The Biot--Savart law reproduces
the velocity field from the vorticity according to the formula:
\begin{equation}
u(x,t)=\frac{1}{4\pi}\int \frac{ (x-y)\wedge \omega(y)}{|x-y|^3}dy\,.
\end{equation}
For the essential part of this kernel, we introduce two smooth
nonnegative radial functions $\beta^\delta_1$ and $\beta^\delta_2$
with
\begin{equation} \beta^\delta_1+\beta^\delta_2=1\,, \beta^1_\delta
=0 \hbox{ for }|x|>2\delta\,\, \hbox{ and }\beta^2_\delta =0 \hbox{
for }|x|<\delta\,.
\end{equation}
Then we have
\begin{eqnarray}
&&|\tilde S (\omega)|\le| \int (\frac{y}{|y|}\cdot\xi(x))({\mathrm{
Det}}(\frac{y}{|y|},
\xi(x+y),\xi(x))\beta^1_\delta(|y|)|\omega(x+y)|\frac{dy}{|y|^3}|
+\nonumber
\\
&& |\int (\frac{y}{|y|}\cdot\xi(x))({\mathrm {Det}}(\frac{y}{|y|},
\xi(x+y),\xi(x)) \beta^2_\delta
(|y|))|\omega(x+y)|\frac{dy}{|y|^3}|\,.
\end{eqnarray}
For the first term we use the bound
\begin{equation}
({\mathrm {Det}}(\frac{y}{|y|}, \xi(x+y),\xi(x)) \beta^1_\delta
(|y|))\le \frac{ k_2(t)}{4\pi} |y|\, .
\end{equation}
to obtain:
\begin{eqnarray}
&&| \int (\frac{y}{|y|}\cdot\xi(x))({\mathrm{
Det}}(\frac{y}{|y|},
\xi(x+y),\xi(x))\beta^1_\delta(|y|)|\omega(x+y)|\frac{dy}{|y|^3}|
\nonumber\\
&&\le k_2(t)\delta |\!|\omega|\!|_\infty
\end{eqnarray}

Next, we write  the second term as
\[\int (\frac{y}{|y|}\cdot\xi(x))({\mathrm {Det}}(\frac{y}{|y|}, \xi(x+y),\xi(x))
\beta^2_\delta (|y|)) (\xi(x+y)\cdot (\nabla_y \wedge
u(x+y)))\frac{dy}{|y|^3}
\]
and integrate by parts with respect to $y$. With the
Lipschitz regularity of $\xi$ one has
\begin{equation*}
|\nabla_y\Bigg((\frac{y}{|y|}\cdot\xi(x))({\mathrm {Det}}(\frac{y}{|y|}, \xi(x+y),\xi(x))\Bigg)|\le C k_2(t)\,.
\end{equation*}
Therefore, one has (observing that the  terms coming from large values of $|y|$  and the terms coming from the derivatives of $\beta^2_\delta(|y|)$    give more regular contributions.)
 \begin{eqnarray*}
&& |\int (\frac{y}{|y|}\cdot\xi(x))({\mathrm {Det}}(\frac{y}{|y|},
\xi(x+y),\xi(x)) \beta^2_\delta
(|y|))|\omega(x+y)|\frac{dy}{|y|^3}|\nonumber\\
&&\le \int |\nabla_y(\frac{y}{|y|}\cdot\xi(x))({\mathrm {Det}}(\frac{y}{|y|},
\xi(x+y),\xi(x)))|\beta^2_\delta
(|y|)) \frac{dy}{|y|^3} |\!| u|\!|_\infty\nonumber\\
&&\le C k_2(t)|\log( \delta) | |\!|u|\!|_\infty\nonumber\\
&&\le  Ck_1(t)k_2(t)| \log \delta| \,.
\end{eqnarray*}
Finally, inserting the above estimates in (\ref{gron}) one
obtains for
$|\!|\omega|\!|_\infty>1$ and $\delta=|\!|\omega|\!|_\infty^{-1}$
\begin{equation*}
\frac{d|\!|\omega|\!|_{\infty}}{dt}\le
Ck_2(t)(1+k_1(t))|\!|\omega|\!|_{\infty} \log|\!|\omega|\!|_\infty\,,
\end{equation*}
and the conclusion follows as in the case of the
Beale--Kato--Majda Theorem.

The reader is referred, for instance, to the book of Majda and Bertozzi
\cite{Majda-Bertozzi} and the recent review of
Constantin \cite{Constantin} for addition relevant material.

\subsection{About the two-dimensional case}

In $2d$ case the vorticity $\omega =\nabla \wedge u $ obeys the
equation
\begin{equation}
\del_t (\nabla \wedge u) +(u\cdot\nabla) (\nabla \wedge u)
=0\,.\label{transvort}
\end{equation}
This evolution equation guarantees the persistence of any  $L^p$
norm ($1\le p\le \infty$) of the vorticity. Taking advantage of this
observation  Youdovitch proved in his remarkable paper \cite{YU} the
existence, uniqueness, and global regularity for all solutions with
initial vorticity in $L^\infty\,.$ If the vorticity is in $L^p$, for
$1<p\le \infty$, then one can prove the existence of weak solutions. The
same results hold also for $p=1$ and for vorticity being a finite
measure with ``simple" changes of sign. The proof is more delicate
in this limit case, cf. Delort \cite{DE} and the section  \ref{ke}
below.

\section{Pathological behavior of solutions }

Continuing with the comments of the previous section one should
recall the following facts.

\noindent$\bullet$ First, in the  three-dimensional case.

i) There is no result concerning the global, in time existence of
smooth solution. More precisely, it is not known whether the
solution of the Euler dynamics defined with initial velocity, say in
$H^s\,,$ for $s>\frac{3}{2} +1$, on a finite time interval can be
extended as a regular, or even as a weak, solution for all positive
time.

ii) There is no result concerning the existence, even for
a small time, of a weak solution for initial data less
regular than in the above case.

iii) Due to the scaling property of the Euler equations in $\R^3$,
the problem of global, in time existence, for small initial data, is
equivalent to the global existence for all initial data and for all
$t\in \R\,.$

\noindent $\bullet$ Second, both in the $2d$ and the $3d$ cases,  the fact that a
function $u\in L^2([0,T]; L^2(\R^d))$ is a weak solution, i.e.,
that it
satisfies the following relations in the sense of distributions
\begin{equation}
\del_t u +\nabla \cdot (u\otimes u)+\nabla p=0\,,\,\,\nabla\cdot u=0, u(x,0)=u_0(x)\,, \label{weak}
\end{equation}
is not enough to define it uniquely in terms of the initial data
(except in $2d$ with the additional regularity assumption that
$\nabla \wedge u_0\in L^\infty\,.)$
 More precisely, one can construct,  following
 Scheffer \cite{SC} and Shnirelman \cite{Sh}, both in  $2d$ and
 $3d$, nontrivial
 solutions $u\in L^2(\R_t;L^2(\R^d))$ of (\ref{weak}) that
 are of compact support in space and time.

  The following examples may contribute to the understanding of the underlying difficulties.
First, one can exhibit (cf. Constantin \cite{CO}, Gibbon and Ohkitani
\cite{SG}, and references therein)  blow up for smooth solutions,
with infinite energy, of the $3d$ Euler equations. Such solutions
can be constructed as follows. The solution $u$ is  $(x_1,x_2)$
periodic on a lattice $(\R/ {L \ZZ})^2$ and is defined for all
$x_3\in \R$ according to the formula
\begin{equation*}
u=(u_1(x_1,x_2,t),u_2(x_1,x_2,t),x_3\gamma(x_1,x_2,t))=(\tilde u ,
x_3 \gamma)\,,
\end{equation*}
which is determined by the following equations:

To maintain the divergence free condition, it is required that
\begin{equation*}
\nabla\cdot \tilde u + \gamma=0\; ,
\end{equation*}
and to  enforce the Euler dynamics, it is required that
\begin{eqnarray*}
&&\del_t (\nabla \wedge \tilde u) +(\tilde u \cdot \nabla) (\nabla \wedge \tilde u )=\gamma  \tilde u \\
&&\del_t \gamma +(\tilde u \cdot \nabla) \gamma = -\gamma^2 + I(t)
\,,
\end{eqnarray*}
and finally to enforce the $(x_1,x_2)$ periodicity it is required
that
\begin{equation*}
I(t)= -\frac 2 {L^2} \int_{[0,L]^2} (\gamma(x_1,x_2,t))^2dx_1dx_2 \,.
\end{equation*}
Therefore, the scalar function $\gamma$ satisfies an
integrodifferential Ricatti equation of the following form
\begin{equation*}
\del_t \gamma +\tilde u \nabla \gamma = -\gamma^2 -\frac 2 {L^2}
\int_{[0,L]^2} (\gamma(x_1,x_2,t))^2dx_1dx_2\,,
\end{equation*}
from which the proof of the blow up, including explicit nature of
this blow up, follows.

The above example can be considered as non-physical
because the initial energy
$$
\int_{(\R^2/L)^2\times\R} |u(x_1,x_2,x_3,0)|^2dx_1dx_2dx_3
$$
is infinite. On the other hand,  it  is instructive because it shows
that the conservation of energy, in the Euler equations, may play a
crucial role in the absence of  singularity. Furthermore, an
approximation of the above solution, by a family of finite energy
solutions, would probably be possible but to the best of our
knowledge this has not yet been done. Such an approximation
procedure  would lead to the idea that no uniform bound can be
obtained for the stability or regularity of $3d$ Euler equations.
Along these lines, one has the following proposition.
\begin{proposition}
For $1<p<\infty$ there is no continuous
function $\tau \mapsto  \phi(\tau)$ such that for
any smooth solution of the Euler equations
the following estimate
\begin{equation*}
|\!| u(\cdot,t)|\!|_{W^{1,p}(\Omega)}\le \phi(|\!|
u(\cdot,0)|\!|_{W^{1,p}(\Omega)})\,,
\end{equation*}
is true.
\end{proposition}
Observe that the above statement is not in contradiction with the
local stability results, which produce local control of higher norm
at time $t$ in term of higher norm at time $0$ as  done in
(\ref{scalarprod}) according to the formula
\begin{equation*}
{\mathrm {for }} \, s>\frac 5 2\,,\quad  |\!|
u(t)|\!|_{H^s(\Omega)}\le \frac { |\!| u(0)|\!|_{H^s(\Omega)}}{1-Ct
|\!| u(0)|\!|_{H^s(\Omega)}}\,.
\end{equation*}
{\bf Proof.}  The proof  is done by inspection of a pressureless
solution,  defined on a periodic box $(\R/\ZZ)^3$ of the form
\begin{equation*}
u(x,t)= (u_1(x_2), 0, u_3(x_1-t u_1(x_2), x_2))\\,
\end{equation*}
which satisfies
\begin{equation*}
\nabla \cdot u =0\,,\,\,\,\, \del_t u +u\cdot \nabla u=0\,.
\end{equation*}
Therefore, the initial data satisfies the relation
\begin{eqnarray}
&&|\!| u(\cdot,0)|\!|^p_{W^{1,p}(\Omega)}\simeq\int_0^1 |\del_{x_2} u_1(x_2)|^pdx_1+\nonumber\\
&&\int_0^1\!\!\int_0^1( |\del_{x_1} u_3(x_1,x_2)|^pdx_1dx_2 +
|\del_{x_2} u_3(x_1,x_2)|^p)dx_1dx_2 \,.\label{finite}
\end{eqnarray}
And for  $t>0$
\begin{eqnarray}
&&|\!| u(\cdot,t)|\!|^p_{W^{1,p}(\Omega)}\simeq\int |\del_{x_2} u_1(x_2)|^pdx_1dx_2dx_3 +\nonumber\\
&&\int_0^1\!\!\int_0^1( |\del_{x_1} u_3(x_1,x_2)|^pdx_1dx_2
+  |\del_{x_2} u_3(x_1,x_2)|^p)dx_1dx_2\nonumber\\
&&+t^p\int_0^1\!\!\int_0^1|\del_{x_2} u_1(x_2)|^p |\del_{x_1}
u_3(x_1,x_2)|^pdx_1dx_2 \,.
 \label{infinite}
 \end{eqnarray}
Then a convenient choice of $u_1$ and  $u_3$ makes the left-hand
side of (\ref{finite}) bounded and the term
$$
t^p\int_0^1\!\!\int_0^1|\del_{x_2} u_1(x_2)|^p |\del_{x_1}
u_3(x_1,x_2)|^pdx_1dx_2\,,
$$
on the right-hand side of (\ref{infinite}) grows to infinity as $t
\to \infty$. The proof is then completed by a regularization
argument.
\begin{remark}
With smooth initial data, the above construction gives an example of
global, in time, smooth solution with vorticity growing (here only
linearly) for $t\rightarrow \infty$.
\end{remark}
As in the case of the Riccati differential inequality $ y'\le C y^2$,
one can obtain sufficient conditions for the existence of a smooth
solution during a finite interval of  time (say $0\le t <T$). On the
other hand, this gives no indication on the possible appearance of
blow up after such time. Complicated phenomena that appear in the
fluid, due to strong nonlinearities, may later interact in such a
way that they balance each other and bring back the fluid to a
smooth regime. Such phenomena is called {\it singularity depletion}.

An example which seems to illustrate such cancelation has been
constructed by Hou and Li \cite{HL}. It is concerning axi-symmetric
solutions of the $3d$ Euler equations form $rf(z)$, which obviously
possess  infinite energy.

Specifically, let us start with the following system of
integro-differential equations with solutions that are defined for
$(z,t)\in (\R/{\mathbf Z})\times \R^+$
\begin{eqnarray}
&&u_t+ 2 \psi u_z=-2vu\,,\,\, v_t+ 2 \psi v_z= u^2-v^2+c(t)\label{preperio}\\
&& \psi_z= v\,,\,\, \int_0^1 v(z,t)dz=0\,.\label{perio}
\end{eqnarray}
In (\ref{preperio}), the $z$ independent function $c(t)$ is chosen to
enforce the second relation of (\ref{perio}), which in turn makes
the function $\psi(z,t)$ $1-$periodic in the $z$ direction. As a
result one has the following:
\begin{lemma}
For any initial data $(u(z,0),v(z,0))\in C^m(\R/ {\mathbf Z})$\,,
with $m\ge 1$, the system (\ref{preperio}) and (\ref{perio}) has a
unique global, in time, smooth solution.
\end{lemma}
{\bf Proof.}  The proof relies on a global {\it a priori}
 estimate. Taking
the derivative with respect to $z$ variable gives (using the
notation $(u_z,v_z) =(u',v')$) :
\begin{eqnarray*}
&&u'_t+ 2 \psi u'_z-2\psi_z u'=-2v'u-2vu'\\
&&v'_t+ 2 \psi v'_z-2\psi_z v'= 2uu'-2vv'\,.
\end{eqnarray*}
Next, one uses the relation $\psi_z=-v$, multiplies the first
equation by $u$, multiplies  the second equation by $v$,
 and adds them  to
obtain
\begin{equation}
\frac12(u_z^2+v_z^2)_t+\psi(u_z^2+v_z^2)_z=0\,. \label{maxprinc}
\end{equation}
The relation (\ref{maxprinc}) provides a uniform $L^\infty $ bound
on the $z-$derivatives of $u$ and $v$. A uniform   $L^\infty$ bound
for $v$ follows from the Poincar\'e inequality, and   finally one
uses for $u$ the following Gronwall estimate
\begin{equation}
|\!|u(z,t)|\!|_{L^\infty}\le |\!|u(z,0)|\!|_{L^\infty}e^{t|\!| (u(z,0))^2+ (v(z,0))^2|\!|_{L^\infty}} \label{gron2}\,.
\end{equation}
\begin{remark} The global existence for solution of the
system  (\ref{preperio}) and (\ref{perio}), with no restriction on
the size of the initial data, is a result of delicate
balance/cancelation, which depends on the coefficients of the
system. Any modification of these coefficients may lead to a blow up
in a finite time of the solutions to the modified system. On the
other hand, the solutions of the system
(\ref{preperio})--(\ref{perio}) may grow exponentially in time.
Numerical simulations performed  by \cite{HL} indicate that the
exponential growth rate in (\ref{gron2}) may get saturated.
\end{remark}
The special structure of  the system (\ref{preperio})--(\ref{perio})
is related to {\it  the $3d$ axi-symmetric Euler equations with
swirl} as follows. Introduce the orthogonal basis
\[e_r=(\frac xr,\frac yr, 0)\,, \,\, e_\theta=(-\frac yr ,\frac xr,0)\,,\,\,e_z=(0,0,1)\,,
\]
and with the solution of the system   (\ref{preperio}) and (\ref{perio}),
construct solutions of the $3d$ ($2+1/2$) Euler equation
according to the following proposition.
\begin{prop}
Assume that $u(z,t) $ and $\psi(z,t)$ are solutions of the systems
(\ref{preperio}) and (\ref{perio}), then the function
\begin{equation*}
U(z,t)=-r\frac{\del \psi (z,t)}{\del z} e_r+ r u(z,t) e_\theta +2
r\psi(z,t)e_z
\end{equation*}
is a smooth solution of the $3d$ Euler, but of an infinite energy.
Moreover, this solution is defined for all time and  without any
smallness assumption on the size of the initial data.
\end{prop}

\section{Weak limit of solutions of the Navier--Stokes Dynamics}

As we have already remarked in the introduction,  for both practical
problems, as well as  for mathematical analysis, it is feasible to
consider the Euler dynamics as the limit of the Navier--Stokes
dynamics, when the viscosity tends to zero. Therefore, this section
is devoted to the analysis of the {\it weak limit,} as
$\nu\rightarrow 0\,,$ of Leray--Hopf type solutions of the
Navier--Stokes equations in $2d$ and $3d$. We will consider only
convergence over finite intervals of time $0<t<T <\infty$. We also
recall that $\nu$ denotes the dimensionless viscosity, i.e., the
inverse of the Reynolds number.

\subsection{Reynolds stresses tensor and dissipative solutions}

As above, we denote by $\Omega$ an open set in $\R^d\,.$ For any
initial data, $u_\nu(x,0)=u_0(x)\in L^2(\Omega)$, and any given
viscosity, $\nu>0$, the pioneer works of Leray  \cite{LE} and Hopf
 \cite{HO} (see also Ladyzhenskaya \cite{LADY03} for a
 detailed survey), which were  later generalized by
 Scheffer \cite{SC}, and by Caffarelli,
Kohn and Nirenberg \cite{CKN},  showed the existence of functions
$u_\nu$ and  $ p_\nu$ with the following properties
\begin{equation}
u_\nu \in L^\infty((0,T);L^2(\Omega))\cap
L^2((0,T);H^1_0(\Omega)),\,\, \hbox{ for every }\,   T \in
(0,\infty). \label{spaces}
\end{equation}
In addition, they satisfy  the Navier--Stokes equations
\begin{eqnarray}
&&\del_t u_\nu + \nabla \cdot (u_\nu \otimes u_\nu)-\nu \Delta u_\nu +\nabla p_\nu=0\,,\\
&&\nabla\cdot u_\nu =0\,, \,\, u_\nu=0 \hbox{ on } \del \Omega\,,
\label{nunsdyn}
\end{eqnarray}
in the sense of distributions. Moreover,  such solutions satisfy the
``pointwise''  energy inequality
\begin{eqnarray}
&&\frac12 \del_t |u_\nu(x,t)|^2 + \nu|\nabla u_\nu(x,t)|^2 + \nonumber \\
&& \nabla\cdot ((u_\nu \otimes u_\nu)(x,t)-\nu
\nabla\frac{|u_\nu(x,t)|^2}2)+ \nabla \cdot
(p_\nu(x,t)u_\nu(x,t))\le 0\, \label{suitsol}
\end{eqnarray}
or in integrated form
\begin{equation}
\frac12 \del_t\int_\Omega |u_\nu(x,t)|^2 dx+\nu\int_\Omega |\nabla
u_\nu(x,t)|^2dx\le 0\,.\label{glob}
\end{equation}

A pair $\{u_\nu, p_\nu\}$ which satisfies
(\ref{spaces}),(\ref{nunsdyn}) and (\ref{suitsol}) is called a
suitable weak solution of the Navier--Stokes equations, in the sense
of Caffarelli--Kohn--Nirenberg. If  it satisfies, however, the
integrated version of the energy inequality (\ref{glob}) instead of
the pointwise energy inequality (\ref{suitsol}) it will then be
called a Leray--Hopf weak solution of the Navier--Stokes equations.

In two-dimensions (or in any dimension but with stronger hypothesis
on the smallness of the size of the initial data with respect to the
viscosity) these solutions are shown to be smooth, unique, and depend
continuously on the initial data. Furthermore, in this case, one has
equality in the relations (\ref{suitsol}) and (\ref{glob}) instead
of inequality.

Therefore, as a result of the above, and in particular the energy
inequality (\ref{glob}), one concludes that, modulo the extraction
of a subsequence, the sequence $\{u_\nu\}$ converges in the
weak$-*$ topology of $L^\infty(\R_t^+, L^2(\Omega))$ to a limit
$\overline{u}$; and the sequence $\{\nabla p_\nu\}$ converges to a
distribution $\nabla \overline{p}$, as $\nu \to 0$, for which the
following holds
\begin{eqnarray}
&&\overline u \in L^\infty(\R_t^+, L^2(\Omega)),
\nabla \cdot \overline u = 0 \hbox{ in }  \Omega\,,\,\, \overline  u\cdot \vec{n}=0 \hbox{ on } \del \Omega\,,\nonumber\\
&& \int_{\Omega}|\overline u(x,t)|^2dx +2\nu\int_0^t\!\!\int_\Omega |\nabla u|^2dxdt\le \int_{\Omega}|\overline u_0(x)|^2dx\,,\nonumber\\
&&\lim_{\nu \rightarrow 0} (u_\nu\otimes u_\nu)= \overline u\otimes
\overline u +
\lim_{\nu \rightarrow 0} ((u_\nu-\overline u)\otimes (u_\nu-\overline u))\,,\label{aubin}\\
&&\del_t \overline{u} + \nabla\cdot  (\overline{u}  \otimes \overline{u} )
+\lim_{\nu \rightarrow 0} \nabla \cdot \Bigg(\overline{u} -u_\nu)\otimes (\overline{u} -u_\nu) \Bigg)+\nabla \overline p=0\label{rtturb}\,.
\end{eqnarray}
Observe that the term
\begin{equation}
RT(x,t)=\lim_{\nu \rightarrow 0}  (\overline{u}(x,t)
-u_\nu(x,t))\otimes (\overline{u}(x,t) -u_\nu(x,t))  \label{detrt}
\end{equation}
is a positive, symmetric, measure-valued tensor. In analogy with (see
below) the statistical theory of turbulence,  this tensor may carry
the name of Reynolds stresses tensor or turbulence tensor. In
particular, certain turbulent regions will correspond to the support
of this tensor.

This approach leads to the following questions.

1. What are the basic properties (if any) of the tensor $RT(x,t)\,?$

2. When does the tensor $RT(x,t)$ identically equal zero?  Or,  what
is equivalent, when does the limit pair
$\{{\overline u}, {\overline p}\}$ satisfies the Euler equations?

3. When does the energy dissipation
$$
\nu \int_0^T\!\!\int_{\Omega} |\nabla u_\nu (x,t)|^2dx dt
$$
tend to zero as $\nu \to 0$?

4. Assuming that $\{{\overline u}, {\overline p}\}$  is a solution
of the Euler equations is such a solution regular enough to imply
the conservation energy?

Hereafter, we will use the following notation for the $L^2-$norm
\[
|\Phi| = \left (\int_\Omega|\Phi(x)|^2 dx\right)^{1/2}.
\]

\begin{remark}
The tensor $RT(x,t)$  is generated by the high frequency spatial
oscillations of the solution. This feature will be explained in more
details in section \ref{ddsp}. Therefore, such behavior should be
intrinsic and, in particular, independent of orthogonal (rotation)
change of coordinates. For instance, in the $2d$ case,
assuming  that
the function $\overline u$ is regular,  the invariance under
rotation implies the relation
\begin{equation*}
RT(x,t)=\alpha (x,t)Id+ \frac12\beta(x,t)(\nabla\overline u + (\nabla \overline u)^T)\,,
\end{equation*}
where $\alpha (x,t)$ and $\beta(x,t)$ are some scalar valued
(unknown) functions. Thus, the equation (\ref{rtturb})  becomes
\begin{equation}
\del_t \overline{u} + \nabla \cdot (\overline{u}  \otimes \overline{u} )
+  \nabla\cdot (\beta(x,t)\frac12(\nabla\overline u + (\nabla \overline u)^T)) +\nabla (\overline p+\alpha(x,t))=0\label{keps}\,.
\end{equation}
Of course, this ``soft information" does not indicate whether
$\beta(x,t)$ is zero or not. It also does not indicate whether this
coefficient is positive,  nor  how to compute it. But this turns out
to be the turbulent eddy diffusion coefficient that is present in
classical engineering turbulence models like the Smagorinsky or the
$k\epsilon$ models (see, e.g., \cite{LASP}, \cite{MOPI},
\cite{Pope},  and \cite{SMA}.
 \end{remark}
\begin{remark}
Assume that the limit $\{{\overline u},  {\overline p}\}$ is a
solution of the Euler equations which is regular enough to ensure
the conservation of energy, i.e. $|{\overline u}(t)|^2=|u_0|^2$.
Then by virtue of the energy relation (\ref{glob}) we have
\begin{equation*}
\frac12|u_\nu(t)|^2+\nu\int_0^t|\nabla u_\nu(s)|^2ds\le
\frac12|u_0|^2\,,
\end{equation*}
and by the weak limit relation
\begin{equation*}
\liminf_{\nu\rightarrow 0} \frac12|u_\nu(t)|^2\ge \frac12|\overline
u(t)|^2\,,
\end{equation*}
one has that the strong convergence and the relation
\begin{equation*}
\liminf_{\nu\rightarrow 0} \nu\int_0^t|\nabla u_\nu(s)|^2ds=0\,
\end{equation*}
hold. The following question was then raised by Onsager \cite{ON}:
``What is the minimal regularity needed to be satisfied by the
solutions of the $2d$ or $3d$ Euler equations that would imply
conservation of energy?". The question was pursued by several
authors up to the contribution of  Eynik \cite{Eynik}, and
Constantin, E and Titi \cite{CET}. Basically in $3d$ it is shown
that if  $u$ is bounded in $L^\infty(\R_t^+,H^\beta(\Omega))\,,$
with $\beta>1/3\,,$ the energy
$$
\frac12\int_\Omega |u(x,t)|^2dx
$$ is constant. On the other hand, arguments borrowed from
statistical theory of turbulence (cf. section \ref{stsp}), show that
the sequence $u_\nu$ will be, in general, bounded in
$L^\infty(\R_t^+,H^{\frac 13}(\Omega))$ and one should observe that
such a statement does not contradict the possibility of decay of
energy in the limit as $\nu\rightarrow 0\,.$

\end{remark}
To study the weak limit of Leray--Hopf  solutions  of the
Navier-Stokes dynamics, P.L.~Lions and R.~Di Perna \cite{DPL}
introduced the notion of {\it Dissipative Solution} of the Euler
equations. To motivate this notion, let $w(x,t)$  be a divergence
free test function, which satisfies $w\cdot\vec n=0$ on the boundary
$\del\Omega\,.$ Let
\begin{equation}
E(w)=\del_t w +P(w\cdot \nabla w)\, ,
\end{equation}
 where  $P$ is  the
Leray--Helmholtz projector (see, e.g., \cite{CF88}). Then for any smooth, divergence free,
solution of the Euler equations $u(x,t)$  in $\Omega\,,$ which
satisfies the boundary condition $u\cdot \vec n =0$ on $\del
\Omega\,,$ one has:
\begin{eqnarray*}
&\del_t u+\nabla \cdot ( u\otimes u)+\nabla p=0\,, \\
&\del_t w+\nabla \cdot ( w\otimes w)+\nabla q=E(w)\,,\\
& \frac{d|u-w|^2}{dt} +2(S(w)(u-w),(u-w))=2(E(w),u-w)\,,
\end{eqnarray*}
where $S(w)$ denotes, as before, the symmetric tensor
$$S(w)=\frac12(\nabla w +(\nabla w)^T)\,.$$
By integration in time this gives
\begin{eqnarray}
&&|u(t)-w(t)|^2\le e^{\int_0^t 2|\!|S(w)(s)|\!|_\infty ds}|u(0)-w(0)|^2 \nonumber\\
&&+2\int_0^te^{\int_s^t 2|\!|S(w)(\tau)|\!|_\infty d\tau}
(E(w)(s),(u-w)(s))ds \label{diss}\,.
\end{eqnarray}
The above observation leads to the following definition
\begin{defi}
A divergence free vector field
$$
u\in w - C(\R _t; (L^2(\Omega))^d),
$$
which satisfies the boundary condition $u\cdot \vec n=0$ on
$\del\Omega$, is called a {\it dissipative solution} 
 of the Euler equations (\ref{Euler1}), if for any smooth
 divergence free vector field $w\,,$ with  $w\cdot \vec n=0$
 on $\del\Omega\,,$ the inequality (\ref{diss}) holds.
\end{defi}
The following statement is easy to verify, but we mention it here
for the sake of clarity.
\begin{theorem}\label{diss-Euler}

i) Any classical solution $u$ of the Euler equations (\ref{Euler1})   is a dissipative solution.

ii) Every dissipative solution satisfies the energy inequality
relation
\begin{equation}
|u(t)|^2\le |u(0)|^2 \label{nonsasha}\,.
\end{equation}

iii) The dissipative solutions are ``stable with
respect to classical solutions".  More precisely,
if  $w$ is a classical solution and $u$ is a dissipative
solution of the Euler equations, then one has
$$
|u(t)-w(t)|^2\le e^{\int_0^t 2|\!|S(w)(s)|\!|_\infty ds}|u(0)-w(0)|^2\,.
$$
In particular, if there exists a classical solution for specific initial data, then any dissipative solution with the same initial data coincides with it.

iv) In the absence of physical boundaries, i.e. in the case of
periodic boundary conditions or in the whole space $\R^d, d=2,3$, any
weak limit, as $\nu\to 0$, of Leray--Hopf solutions of the
Navier--Stokes equations is a dissipative solution of Euler
equations.

\end{theorem}

{\bf Proof and remarks.} The point i) is a direct consequence of the
construction.  To prove ii) we consider $w\equiv0$ as a classical
solution. As a result, one obtains for any dissipative solution, the
relation (\ref{nonsasha}), which justifies the name dissipative.
Furthermore, it shows that the pathological examples constructed by
Scheffer \cite{SC} and Shnirelman \cite{Sh} are not dissipative
solutions of Euler equations.

For the point iii)  we use in (\ref{diss}) the fact that $w$ being a
classical solution implies that $E(w)\equiv 0.$ We also observe that
this statement is in the spirit of the ``weak with respect to
strong" stability result of Dafermos \cite{DA} for hyperbolic
systems.

Next, we prove iv) in the absence of physical boundaries. Let
 $u_\nu$ be a Leray--Hopf solution  of the Navier--Stokes
 system, which satisfy an energy inequality in (\ref{glob}),
 and let $w$ be a classical solution of the Euler equations.
 By subtracting the following two equations from each other
\begin{eqnarray*}
&&\del_t u_\nu +\nabla\cdot (u_\nu \otimes u_\nu)-\nu \Delta u_\nu +\nabla p_\nu=0\\
&&\del_t w +\nabla\cdot (w\otimes w)-\nu \Delta w +\nabla
p=-\nu\Delta w,
\end{eqnarray*}
and  taking the $L^2$ inner product of the difference with
$(u_\nu-w)$ one obtains
 \begin{eqnarray}
&&\frac{d|u_\nu-w|^2}{dt} +2(S(w)(u_\nu-w),(u_\nu-w))\nonumber-2\nu  (\Delta(u_\nu-w),(u_\nu-w)) \\
&&\le 2(E(w),u_\nu-w)-(\nu\Delta w , (u_\nu-w))\label{lay1}\,.
\end{eqnarray}
We stress that the above step is formal, and only through rigorous
arguments one can see the reason for obtaining an inequality in
(\ref{lay1}), instead of an equality. However, this should not be a
surprise because we are dealing with Leray--Hopf solutions, $u_\nu$,
of the Navier--Stokes system which  satisfy an energy inequality in
(\ref{glob}), instead of an equality.

Now, to conclude our proof we observe that in the absence of
physical boundaries one uses the relation
\begin{equation}
-\nu \int \Delta(w-u_\nu)(x,t)\cdot(w-u_\nu)(x,t)dx=\nu \int
|\nabla(w-u_\nu)(x,t)|^2dx
\end{equation} \label{lay230}
and the result follows by letting $\nu$ tend to zero.

\begin{remark} The above theorem states in particular,
and in the absence of physical boundaries, that as long as a smooth
solution of the Euler equations does exist, it is the limit, as $\nu
\to 0$,  of any sequence of Leray--Hopf solutions of the
Navier--Stokes equations with the same initial data. In a series of
papers, starting with Bardos, Golse and Levermore \cite{BGL},
connections between the notion of Leray--Hopf solutions for the
Navier--Stokes equations and {\it renormalized solutions of the
Boltzmann equations,} as defined by P.L.~Lions and Di Perna, were
established. In particular, it was ultimately shown by Golse and
Saint Raymond \cite{GSR} that, modulo the extraction of a
subsequence, and under a convenient space time scaling, any sequence
of such renormalized solutions of the Boltzmann equations converge
(in some weak sense) to a Leray--Hopf solution of the Navier--Stokes
system. On the other hand, it was shown by Saint Raymond \cite{SR}
that, under a scaling which reinforces the nonlinear effect
(corresponding at the macroscopic level to a Reynolds number going
to infinity), any sequence (modulo extraction of a subsequence) of
the renormalized solutions of the Boltzmann equations converges to a
dissipative solution of the Euler equations. Therefore, such a
sequence of normalized solutions of the Boltzmann equations
converges to the classical solution of Euler equations, as long as
such solution exists. In this situation, one should observe that,
with the notion of dissipative solutions of Euler equations,
classical solutions of the Euler equations play a similar role for the
``Leray--Hopf limit" and the ``Boltzmann limit."
\end{remark}

\begin{remark}
There are at least two situations where the notion of dissipative
solution of Euler equations is not helpful.

The first situation is concerned with the $2d$ Euler equations. Let
$u_\epsilon(x,t)$ be the sequence of solutions of the $2d$ Euler
equations corresponding to the sequence of smooth  initial data
 $u_\epsilon(x,0)$. Suppose that the sequence of initial
 data $u_\epsilon(x,0)$
converges weakly, but not strongly, in $L^2(\Omega)$ to an initial
data $\overline u(x,0)$, as $\epsilon\rightarrow 0\,.$ Then for any
smooth, divergence free vector field $w$, one has, thanks to
(\ref{diss}), the relation
\begin{eqnarray}
&&|u_\epsilon (t)-w(t)|^2\le e^{\int_0^t 2|\!|S(w)(s) |\!|_\infty ds}|u_\epsilon (0)-w(0)|^2 \nonumber\\
&&+2\int_0^te^{\int_s^t 2|\!|S(w)(\tau)|\!|_\infty d\tau}
(E(w),u_\epsilon-w)(s)ds \label{diss2}\,.
\end{eqnarray}
However, with the weak convergence as $\epsilon\rightarrow 0$,
 one has only
\begin{equation*}
|\overline{ u  } (0)-w(0)|^2 \le \liminf_{\epsilon\rightarrow
0}|u_\epsilon (0)-w(0)|^2\,,
\end{equation*}
and  (\ref{diss2}) might not hold at the limit, as
$\epsilon\rightarrow 0$.  To illustrate this situation, we consider a
sequence of oscillating solutions of the $2d$ Euler equations of the
form
\begin{equation*}
u_\epsilon(x,t)=U(x,t,\frac{\phi(x,t)}\epsilon) +O(\epsilon)\,,
\end{equation*}
where the map $\theta\rightarrow U(x,t,\theta ) $ is a nontrivial $1-$periodic function. In Cheverry \cite{CHE} a specific example was constructed such that
$$
\overline u = w-\lim_{\epsilon \rightarrow 0} u_\epsilon=
\int_0^1U(x,t,\theta)d\theta
$$
is no longer a solution of the Euler equations. The obvious reason
for that (in comparison to the notion of dissipative solution),
 is the fact that
$$
U(x,0,\frac{\phi(x,0)}\epsilon)
$$
does not converge strongly in $L^2(\Omega)\,.$

 The second situation, which will be discussed at length below,
 corresponds to the weak limit of solutions of the
 Navier--Stokes equations in a domain with physical
 boundary,  subject to the no-slip Dirichlet boundary condition.

As we have already indicated in Theorem \ref{diss-Euler}, one of the
most important features of the above definition of dissipative
solution of Euler equations is that it coincides with the classical
solution of Euler equations, when the latter exists. This can be
accomplished by replacing $w$ in (\ref{diss}) with this classical
solution of Euler equations. Therefore, any procedure for
approximating dissipative solutions of Euler must lead to, in the
limit, to  inequality (\ref{diss}). Indeed, in the absence of
physical boundaries,  we have been successful in showing, in an
almost straight forward manner, in Theorem \ref{diss-Euler}, that the
Leray--Hopf weak solutions of the Navier--Stokes equations converge
to dissipative solutions of Euler equations. On the other hand, in
the presence of physical boundaries, the proof does not carry on in
a smooth manner because of the boundary effects. Specifically, in the
 case of domains with physical boundaries,
 inequality (\ref{lay1}) leads to
 \begin{eqnarray}
&&\frac12 \frac{d|u_\nu-w|^2}{dt} +(S(w)(u_\nu-w),(u_\nu-w))+\nu \int |\nabla(w-u_\nu)|^2dx\nonumber\\
&& \le  (E(w),u_\nu-w)-\nu(\Delta w,  (u_\nu-w))+ \nu \int_{\del
\Omega } \del_n{u_\nu} \cdot w d\sigma\,.\label{lay23}
\end{eqnarray}
The very last term in (\ref{lay23}) represents the boundary effect.
We will discuss below the subtleties in handling this term.

\end{remark}


\section{No-slip Dirichlet boundary conditions for  the Navier--Stokes dynamics}

This section is devoted to the  very few available results
concerning the limit, as  $\nu \rightarrow 0$, of solutions of the
Navier--Stokes equations in a domain $\Omega \subset \R^d \, , \quad
d=2, 3$ with the homogenous (no-slip)  Dirichlet boundary condition
$u_\nu=0$ on $\del \Omega$. This boundary condition  is not the
easiest to deal with, as far as the  zero viscosity limit is
concerned. For instance,  the solutions of the of $2d$
Navier--Stokes equations, subject to  the boundary conditions
$u_\nu\cdot n =0 $ and $\nabla \wedge u_\nu =0$, are much better
understood and much easier to analyze mathematically \cite{BA} as
the viscosity $\nu \to 0$. However, the no-slip boundary condition
is the one which is more suitable to consider physically for the
following reasons.

i)  It can be deduced in the smooth (laminar) regime, from the
Boltzmann kinetic equations when the interaction with the boundary
is described by a scattering kernel.

ii)  It generates the pathology that is observed in physical
experiments, like the  Von Karman vortex streets. Moreover, one
should keep in mind that almost all high Reynolds number  turbulence
experiments involve a physical boundary (very often turbulence is
generated by a pressure driven flow through a grid!)

The problem emerges first from the boundary layer. This is because
for the Navier--Stokes dynamics, the whole velocity field
equals zero
on the boundary, i.e. $u_\nu=0$ on $\del \Omega$,
while for the Euler dynamics  it is only the normal
component of velocity field is equal to zero on the boundary, i.e. $u
\cdot {\vec n}=0$ on  $\del \Omega$. Therefore, in the limit, as the
viscosity $\nu \to 0$, the tangential component of the velocity
field of the Navier--Stokes dynamics,  $u_\nu$, generates, by its
``jump'', a boundary layer. Then, unlike the situation  with linear
singular perturbation problems, the nonlinear advection term of the
Navier--Stokes equations  may propagate this instability inside the
domain.

As we have already pointed out, the very last term in (\ref{lay23}),
i.e. the boundary integral term in the case of no-slip boundary
condition,
\begin{eqnarray}
&&\nu \int_{\del\Omega} \frac{\del u_\nu}{\del n}\cdot w d\sigma=
\nu \int_{\del\Omega} (\frac{\del u_\nu}{\del n})^\tau \cdot w^\tau d\sigma\nonumber\\
&&\hskip 1in = \nu \int_{\del\Omega} (\nabla \wedge u_\nu) \cdot
(\vec n\wedge w) d\sigma\,, \label{lay4}
\end{eqnarray}
is possibly responsible for the loss of regularity in the
limit as
$\nu\to 0$. This is stated more precisely in the following.
\begin{proposition}
Let $u(x,t)$ be a solution of the incompressible Euler equations in
$\Omega \times (0,T]\,,$ with the following regularity assumptions.
\begin{eqnarray*}
&&S(u)=\frac12(\nabla u +(\nabla u)^T)\in
L^1((0,T);L^\infty(\Omega)),\,\,
{\hbox {and}}\\
&&u\in L^2((0,T); H^s(\Omega)),\,\, {\hbox {for}} \,\, s>\frac 12\,.
\end{eqnarray*}
Moreover, suppose that the sequence $u_\nu$, of Leray--Hopf
solutions of the Navier--Stokes dynamics (with no-slip boundary
condition) with the initial data $u_\nu(x,0)=u(x,0)$, satisfies the
relation
\begin{equation*}
\lim_{\nu \rightarrow 0} \nu |\!| P_{\del\Omega}(\nabla \wedge
u_\nu)|\!|_{L^2((0,T); H^{-s+\frac12}(\del \Omega))}=0\,,
\end{equation*}
where $P_{\del\Omega}$ denotes the projection on the tangent plane
to $\del\Omega$ according to the formula.
\begin{equation*}
P_{\del\Omega} (\nabla \wedge u_\nu)= \nabla \wedge u_\nu- ((\nabla
\wedge u_\nu) \cdot \vec n)\vec n\,.
\end{equation*}
Then, the sequence $u_\nu$ converges to $u$ in $C((0,T);
L^2(\Omega))$.
\end{proposition}
The proof is a direct consequence of (\ref{lay23}) and (\ref{lay4}),
with $w$ being replaced by $u$. This Proposition  can be improved
with the following simple and beautiful theorem of Kato which takes
into account the  vorticity production in the boundary layer $\{
x\in \Omega \,| d (x,\del \Omega)<\nu\}\,,$ where $d(x,y)$ denotes
the Euclidean distance between the points $x$ and $y$.

\begin{theorem}\label{kato}
Let $u(x,t)\in W^{1, \infty} ((0, T)\times \Omega)$ be a solution of
the Euler dynamics, and let $u_\nu$ be a sequence of Leray--Hopf
solutions  of the Navier--Stokes dynamics with no-slip boundary
condition.
{\begin{eqnarray} \del_t u_\nu -\nu \Delta u_\nu
+\nabla\cdot ( u_\nu \otimes  u_\nu) +\nabla p_\nu=0 \,, \,
u_\nu(x,t)=0\hbox { on } \del \Omega\,,\label{NSB} \label{EUB}
\end{eqnarray}}
with initial data $u_\nu(x,0)=u(x,0)$\,. Then, the following facts
are equivalent.
 \begin{eqnarray}
&&(i) \quad  \lim_{\nu \rightarrow 0} \nu\int_0^T \int_{\del\Omega} (\nabla \wedge u_\nu) \cdot (\vec n\wedge u) d\sigma dt=0\label{ib}\\
&& (ii) \quad u_\nu (t)\rightarrow   u(t) \hbox{ in }  L^2(\Omega)  \hbox{ uniformly in}\, t \in [0,T]\label{iib}\\
&& (iii) \quad u_\nu (t)\rightarrow  u(t)  \hbox{ weakly in } L^2(\Omega) \hbox{ for each } t \in [0,T] \label{iiib}\\
&& (iv) \quad \lim_{\nu \rightarrow 0} \nu \int_0^T\int_{\Omega}
|\nabla u_\nu (x,t) |^2dxdt =0
\label{iiiib}\\
&& (v) \quad\lim_{\nu \rightarrow 0}\nu \int_0^T\int_{\Omega \cap
\{d(x,\del\Omega )<\nu\} }| \nabla u_\nu (x,t) |^2dxdt =
0\,.\label{v}
\end{eqnarray}
\end{theorem}

\noindent{\bf Sketch of the proof.} The statement
 (\ref{iib}) is deduced from (\ref{ib})  by replacing
 $w$ by $u$ in (\ref {lay23}) and (\ref{lay4}). No
 proof is needed to deduce (\ref{iiib}) from (\ref{iib}),
 or (\ref{iiiib}) from (\ref{iiib}).

Next, one recalls the energy inequality (\ref{glob}), satisfied by
the Leray--Hopf solutions of the the Navier--Stokes dynamics
\begin{equation}
\frac12\int_\Omega |u_\nu(x,T)|^2dx +\nu\int_0^T\int_\Omega |\nabla
u_\nu(x,t)|^2dxdt \le
\frac12\int_\Omega|u(x,0)|^2dx\,.\label{encorenerg}
\end{equation}
By virtue of the weak convergence, as stated in (\ref{iiib}), and
the fact $u$ is a smooth solution of the Euler dynamics, one has
\begin{eqnarray}
&&\lim_{\nu\rightarrow 0} \frac12\int_\Omega |u_\nu(x,T)|^2dx\ge
\frac12\int_\Omega | \lim_{\nu\rightarrow 0}
u_\nu(x,T)|^2dx\nonumber\\
&& \hskip 0.5in = \frac12\int_\Omega |u (x,T)|^2dx=
\frac12\int_\Omega |u (x,0)|^2dx\,. \label{regularsol}
\end{eqnarray}
Together with (\ref{encorenerg}) this shows that (\ref{iiib})
implies (\ref{iiiib}).

The most subtle part in the proof of this theorem is the fact that
(\ref{v}) implies (\ref{ib}). The first step is the construction of
a divergence free function $v_\nu(x,t)$ with support in the region
$\{x\in\Omega | d(x, \del \Omega )\le \nu\} \times [0,T)$, which
coincides with $u$ on $\del\Omega \times[0,T]$, and which satisfies
(with $K$ being a constant that is independent  of $\nu$) the
following estimates
\begin{eqnarray}
&&|\!|v_\nu|\!|_{L^\infty(\Omega\times (0,T))}+ |\!|d (x,\del \Omega ) \nabla v_\nu|\!|_{L^\infty( \Omega \times (0,T))}\le K\,,\\
&& |\!|(d(x,\del \Omega ) )^2\nabla v_\nu|\!|_{L^\infty( \Omega \times (0,T))}\le K \nu \,,\\
&&|\!|v_\nu|\!|_{L^\infty((0,T) ;L^2(\Omega) )}+|\!|\del_tv_\nu|\!|_{L^\infty((0,T);L^2(\Omega) )}\le K \nu^{\frac 12}\,,\\
&&|\!|\nabla v_\nu|\!|_{L^\infty((0,T) ;L^2(\Omega) )}\le K \nu^{-\frac 12}\,,\\
&&|\!|\nabla v_\nu|\!|_{L^\infty ( \Omega \times (0,T))}\le
K\nu^{-1}\,.
\end{eqnarray}

Then we multiply the  Navier--Stokes equations by  $v_\nu$ and
integrate to obtain
\begin{eqnarray}
&&-\nu\int_0^T \int_{\del\Omega} (\nabla \wedge u_\nu) \cdot (\vec n\wedge u) d\sigma dt=-\nu\int_0^T \int_{\del\Omega}\frac{\del u_\nu}{\del n} \cdot v_\nu  d\sigma dt \nonumber\\
&&=-\nu \int_0^T \int_\Omega \Delta u_\nu \cdot v_\nu dx dt -\nu\int_0^T \int_\Omega (\nabla u_\nu : \nabla v_\nu) dx dt\\
&&{\hskip -0.5in }=-\int_0^T(\del_t u_\nu, v_\nu)dt -
\int_0^T(\nabla \cdot (u_\nu\otimes u_\nu) , v_\nu)dt  -\nu
\int_0^T(\nabla u_\nu,\nabla v_\nu)dt\,.\label{kato1}
\end{eqnarray}
Eventually, by using  the above estimates and (\ref{v}), one can show
that
\begin{eqnarray}
\lim_{\nu\rightarrow 0} \left(\int_0^T((\del_t u_\nu, v_\nu)
+(\nabla \cdot (u_\nu\otimes u_\nu) , v_\nu)+\nu (\nabla
u_\nu,\nabla v_\nu))dt\right)=0\,,
\end{eqnarray}
which completes the proof.
\begin{remark}
In Constantin and Wu \cite{Constantin-Wu} the authors study the
rate of convergence of solutions of the $2d$ Navier--Stokes equations
to the solutions of the Euler equations in the absence of physical
boundaries, for finite intervals of time.  Their main
observation is that while the rate of convergence, in the $L^2$--norm,
 for smooth initial
data is of the order $\mathcal{O}(\nu)$  it is of the order of
$\mathcal{O}(\sqrt{\nu})$ for less smooth initial. The order of
convergence $\mathcal{O}(\sqrt{\nu})$ is, for instance, attained
for the vortex patch data with smooth boundary. In this case the fluid
develops an internal boundary layer which is responsible for this
reduction in the order of convergence.

In the $2d$ case and for  initial data of finite $W^{1,p}$ norm,
with $p>1$, (and also for initial data with vorticity being a finite
measure with a ``simple" change of sign \cite{DE}, \cite{LNZ}) one can
prove the existence of ``weak  solutions of the Euler dynamics." In
the absence of physical boundaries, such solutions are limit points
of a family of (uniquely determined) Leray--Hopf solutions of the
Navier--Stokes dynamics. However, these weak solutions of the, Euler
equations are not uniquely determined, and the issue of the
conservation of energy for these weak solutions is, to the best of
our knowledge, completely open.

On the other hand, in a domain with physical  boundary and with
smooth initial data,  Theorem \ref{kato}  shows a clear cut
difference between the following two situations (the same remark
being valid locally in time for the $3d$ problems).

i) The mean rate of dissipation of energy
$$
\epsilon= \frac{\nu}{T}\int_0^T \int|\nabla u_\nu(x,t)|^2dxdt
$$
goes to zero as  $\nu\rightarrow 0$, and the sequence $u_\nu$ of
Leray--Hopf solutions converges strongly to the regular solution
$\overline u$ of the Euler dynamics.

ii) The mean rate of  dissipation of energy does not go to zero as
$\nu \to 0$   (modulo the extraction of a subsequence), so the
corresponding weak limit $u$ of $u_\nu$ does not conserve energy,
i.e.
\begin{equation*}
\frac12 |u(t)|^2<\frac12 |u(0)|^2\,, \quad {\text {for some}}\quad t
\in (0,T)\,,
\end{equation*}
and one of the following two scenarios may occur  at the limit.

a) At the limit one obtains a  weak solution
(not a strong solution) of the Euler dynamics that
exhibits energy decay. Such a scenario is compatible with a uniform
estimate for the Fourier spectra
$$
E_\nu(k,t)=  |\hat{ u}_\nu(k,t)|^2 |k|^{d-1}\,,\,\,\hat
{u}_\nu(k,t)=\frac1{(2\pi )^d}\int_{\R^d} e^{-ikx }u_\nu(x,t)dx \,,
$$
which may satisfy a  uniform, in $\nu$, estimate of the following
type
\begin{equation}
E_\nu(k,t)\le C |k|^{-\beta}\,, \label{kolm}
\end{equation}
provided $\beta <5/3$. Otherwise, this would be in contradiction
with the results of  Onsager \cite{ON}, Eyink \cite{Eynik}, and
Constantin, E  and Titi  \cite{CET}.

b) No estimate of the type (\ref{kolm}) is uniformly (in the
viscosity)  true, and the limit is not even a solution of the Euler
dynamics, rather a solution of a  modified system of equations with a
term related to  turbulence modeling - an ``eddy-viscosity'' like
term.
\end{remark}

\section{ Deterministic and statistical spectra of turbulence}

\subsection{Deterministic spectra and Wigner transform} \label{ddsp}

The purpose of this section is the introduction of Wigner measures
for the analysis (in dimensions $d=2,3$) of the Reynolds stresses
tensor
\begin{equation}
RT(u_\nu)(x,t)=\lim_{\nu \rightarrow 0} ((u_\nu-\overline u)\otimes
(u_\nu-\overline u))\,,
\end{equation}
which appears in the weak limit process of solutions of the
Navier--Stokes equations $u_\nu$, as $\nu \to 0$, cf. (\ref{aubin}).
(Notice that $RT(u_\nu)(x,t)$ is independent of the viscosity $\nu$,
but it depends on the sequence $\{u_\nu\}$.) This point of view will
be compared below (cf. section \ref{stsp}) to ideas emerging from
statistical theory of turbulence.

Let $\{u_\nu,p_\nu\}$ be a sequence of  Leray--Hopf solutions of the
Navier--Stokes equations, subject to no-slip Dirichlet boundary
condition in a domain $\Omega$ (with a physical boundary).  Thanks
to the energy inequality (\ref{glob}) (possibly equality in some
cases)

\begin{eqnarray}
\frac 12 |u_\nu(\cdot,t)|^2 +\nu \int_0^t |\nabla
u_\nu(\cdot,t)|^2dt \le \frac 12 |u(\cdot,0)|^2 \,, \label{energy}
\end{eqnarray}
the sequence  $\{u_\nu\}$  converges (modulo
the extraction of a subsequence), as $\nu \to 0$, in the weak$-*$ topology of the
Banach space $L^\infty ((0,T);L^2(\Omega))$, to a divergence free
vector field $\overline u\,,$ and the sequence of distributions
$\{\nabla p_\nu\}$ converges to a distribution $\nabla \overline
p\,.$ Moreover, the pair $\{\overline u, \overline p \}$ satisfies
the system of equations
\begin{eqnarray}
&&\nabla \cdot  \overline{u}=0\, \,\hbox{ in } \Omega\,,\quad \overline u \cdot \vec n=0\, \hbox{ on } \del\Omega\,,\\
&&\del_t \overline u+ \nabla \cdot (\overline u \otimes \overline u)
+ \nabla \cdot RT(u_\nu)+ \nabla \overline p =0 \,
\,\,\hbox{in}\,\,\Omega\label{wig1}\,,
\end{eqnarray}
(cf. (\ref{rtturb})). Being concerned with the behavior of the
solution inside the domain, we consider an arbitrary open subset
$\Omega'$\,, whose closure is a compact subset of $\Omega$, i.e.
$\overline{\Omega'}\subset\!\subset\Omega$. Assuming  the
weak$-*$ limit function $\overline u$ belongs to the space
$L^2((0,T);H^1(\Omega))$ we introduce the function
\begin{equation*}
v_\nu = a(x)(u_\nu -\overline u)\,,
\end{equation*}
with  $a(x)\in {\mathcal D} (\Omega)\,, a(x)\equiv 1$ for all $ x\in
\Omega'\,.$ As a result of (\ref{energy}), and the above
assumptions,  the sequence $v_\nu$ satisfies the uniform estimate
\begin{equation}
\nu \int_0^\infty \!\!\int_\Omega |\nabla v_\nu|^2dx \le C\,.
\end{equation}
Consequently, the sequence $v_\nu$ is (in the sense of Gerard,
Mauser, Markowich and Poupaud \cite{GMMP} ) $\sqrt \nu-$oscillating.
Accordingly, we introduce the deterministic correlation spectra, or
{\it Wigner transform}, at the scale  $\sqrt \nu$\,:
\begin{equation*}
\widehat{ RT(v_\nu)} (x,t,k)= \frac1 {(2\pi) ^d}\int_{\R^d_y}
e^{ik\cdot y} ( v_\nu(x-\frac {\sqrt {\nu}}2 y)\otimes v_\nu(x+\frac
{\sqrt {\nu}}2 y))dy\,.
\end{equation*}
By means of  the inverse Fourier transform, one has
\begin{equation}
v_\nu(x,t)\otimes v_\nu(x,t) =
\int_{\R^d_k}\widehat{ RT(v_\nu)} (x,t,k)dk\,. \label{oscilb}
\end{equation}
The tensor $\widehat{ RT(v_\nu)} (x,t,k)$ is  the main object of
section 1 of \cite{GMMP}. Modulo the extraction of a subsequence,
 the tensor $\widehat{ RT(v_\nu)} (x,t,k)$ converges weakly, as $\nu \to
0$, to a nonnegative symmetric matrix-valued measure $\widehat
{RT}(x,t,dk)$, which is called a {\it Wigner measure or Wigner
spectra}. Moreover, inside the open subset $\Omega'$ the weak  limit
$\overline u$ is a solution of the equation
\begin{equation*}
\del_t  \overline{u}+ \nabla \cdot (\overline{u}  \otimes
\overline{u}) +\nabla \cdot \int_{\R^d_k} \widehat{ RT}
(x,t,dk)+\nabla \overline p=0\,.
\end{equation*}
The Wigner spectra
has the following properties

i) It is defined by a two points correlation formula

ii)  It is an $(x,t)$   locally dependent object.
Specifically, for
any $\phi\in{\mathcal D(\Omega)}$ one has:
\begin{eqnarray}
&&\lim_{\nu \rightarrow 0}\Bigg(\frac1 {(2\pi)^d}\int_{\R^d_y} e^{ik\cdot y}( (\phi v_\nu)(x-\frac {\sqrt {\nu}}2 y)\otimes (\phi v_\nu)(x+\frac {\sqrt {\nu}}2 y)) dy  \Bigg)\nonumber\\
&&=|\phi(\cdot)|^2 \widehat{RT}(x,t,dk) \label{local}
\end{eqnarray}
Therefore, the construction of $\widehat{RT}(x,t,dk)$ is
independent of the choice of the pair $a(x)$ and the open subset
$\Omega'$.

iii) It is a {\it criteria for  turbulence}: Points $(x,t)$ around
which  the sequence $u_\nu$ remains smooth and converges locally
strongly to $\overline u$, as $\nu \to 0$, are characterized by the
relation
\begin{equation}
{\rm Trace} (\widehat{RT}(x,t,dk))=0\,.
\end{equation}
iii) It is a {\it microlocal object.} In fact, it depends only on
the behavior of the Fourier spectra of the sequence $\phi(x) v_\nu$
(or in fact $\phi(x) u_\nu$) in the  {\it frequency band }
$$
A\le |k|\le \frac  B{\sqrt{\nu}}\,.
$$
\begin{prop} For any pair of strictly positive constants $(A,B)$
and any test functions
$(\psi(k),\phi(x),\theta(t))\in C_0^\infty(\R^d_k)\times
C_0^\infty(\R^d_x)\times C_0^\infty(\R^+_t)$ one has
\begin{eqnarray}
&&\int_0^\infty \!\!\int \psi(k)|\phi(x)|^2\theta(t)  {\rm Trace}(\widehat{ RT(u_\nu)})(x,t,dk)dxdt\nonumber\\
&&=\lim_{\nu\rightarrow 0} \int_0^\infty \theta(t) \int_{A\le |k|\le \frac{B}{\sqrt \nu}} ( \psi(\sqrt \nu k)\widehat{(\phi v_\nu)},\widehat{(\phi v_\nu)})dkdt\label{freq2}\,.
\end{eqnarray}
\end{prop}
The only difference between this presentation and what can be
found in \cite{GMMP} comes from the fact that the weak limit
$\overline u$ has been subtracted from the sequence $u_\nu$.
Otherwise the formula (\ref{oscilb}) together with the energy
estimate are the first statement of the Proposition 1.7 of
\cite{GMMP}; while the formula (\ref{freq2}) is deduced from (1.32)
in \cite{GMMP} by observing that the weak convergence of
$u_\nu-\overline u $ to $0$ implies that
\[
\lim_{\nu\rightarrow 0} \left (\int_0^\infty \theta(t) \int_{|k|\le
A}( \psi(\sqrt \nu k)\widehat{(\phi v_\nu)},\widehat{(\phi
v_\nu)})dkdt \right) =0\,.
\]

\subsection{ Energy spectrum in statistical theory of turbulence}\label{stsp}
The Wigner spectra studied in the previous section turns out to be
the deterministic version of the ``turbulent spectra", which is a
classical concept in the statistical theory of turbulence. The two
points of view can be connected with the introduction of homogenous
random variables. Let $(\mathfrak{M},\mathfrak{F}, dm)$ be an
underlying probability space. A random variable $u(x,\mu)$ is said
to be homogenous if for any function $F$ the expectation of
$F(u(x,\mu))$, namely,
$$
\langle F(u(x,\cdot))\rangle= \int_{\mathfrak{M}}F(u(x,\mu)) dm(\mu)
$$  is independent of $x$, that is
\[
\nabla_x( \langle F(u(x,\cdot) \rangle) =0\,.
\]
In particular, if $u(x,\mu)$ is a homogeneous random vector-valued
function, one has
\begin{equation}
\langle u(x+r,\cdot)\otimes u(x,\cdot)\rangle =\langle u(x+\frac {r}
{2},\cdot)\otimes u(x-\frac {r}{ 2},\cdot)\rangle \,,
\end{equation}
which leads to the following proposition.
\begin{proposition}
Let $u(x,\mu)$ be a homogenous random variable
and denote by $\hat
u(k,\mu)$ its Fourier transform. Then one has:
\begin{equation}
\langle \hat u(k,\cdot)\otimes\overline {\hat u(k,\cdot)}\rangle=
\frac{1}{(2 \pi )^d}\int_{\R^d} e^{-ik\cdot r} \langle u(x+\frac r 2
,\cdot)\otimes u(x-\frac r 2 ,\cdot)\rangle dr \,. \label{spec}
\end{equation}
\end{proposition}
{\bf Proof.} The proof will be given for a homogenous random
variable which is periodic with respect to the variable $x$, with
basic periodic box of size $2\pi L$. The formula (\ref{spec}) is
then deduced by letting $L$ go to infinity. From the Fourier series
decomposition in $(\R/ 2 \pi L \ZZ)^d$, one has
\begin{eqnarray*}
&& \hat u (k,\mu)\otimes \overline {\hat u
(k,\mu)}=\\
&&\frac{1}{(2\pi L)^{2d}}\int_{(\R/2\pi L\ZZ)^d}\int_{(\R/2\pi
L\ZZ)^d} (u (y,\mu) e^{-i\frac{k\cdot y}{ L}}\otimes u (y',\mu)
e^{i\frac{k\cdot y'}L})dy'dy=\\
&&\frac{1}{(2\pi L)^{2d}}\int_{(\R/2\pi L\ZZ)^d} e^{-i\frac{k \cdot
r}L} \int_{(\R/2\pi L\ZZ)^d} (u (y,\mu) \otimes u (y+r,\mu))dy dr\,.
\end{eqnarray*}
Averaging with respect to the probability measure $dm$, using the
homogeneity of the random variable $u (x,\mu)$
and integrating with respect to $dy$, gives us the following
\begin{equation}
\langle  \hat u (k, \cdot)\otimes \overline {\hat u (k, \cdot)}
\rangle = \frac{1}{(2\pi L)^d}  \int_{(\R/2\pi L \ZZ)^d} e^{-i\frac
{k \cdot r}{L} } \langle u_\nu(y+\frac r 2, \cdot) \otimes
u_\nu(y-\frac r2, \cdot)\rangle dr\,.\label{wigal}
\end{equation}
This concludes our proof.

 Next,  assuming that in addition to homogeneity, the expectation of
the two points correlation tensor, $\langle u(x+r,\cdot)\otimes
u(x,\cdot)\rangle$, is {\it isotropic} (i.e., it does not depend on
the direction of the vector $r$, but only on its length) one obtains
the following formula
\begin{eqnarray*}
&&\langle  \hat u (k,\cdot)\otimes \overline {\hat u (k,\cdot)}
\rangle = \frac{1}{(2\pi L)^d}  \int_{(\R/2 \pi L \ZZ)^d} e^{-i\frac
{k \cdot r}L }
\langle u_\nu(y+\frac r 2, \cdot) \otimes u_\nu(y-\frac r2, \cdot)\rangle  dr\nonumber\\
&&= \frac{E(|k|)}{S_{d-1}  |k|^{d-1}}(I-\frac{k\otimes k}{|k|^2})\,,
\end{eqnarray*}
with $S_1= 2\pi\,, S_2=4\pi$, which defines the turbulent spectra
$E(|k|)\,.$

 The notion of homogeneity implies that solutions of the
 Navier--Stokes equations
 satisfy a local version of the energy
 balance, often called the Karman--Howarth
 relation cf. (\ref{Karman}). Specifically, let $\{u_\nu, p_\nu\}$
 be  solutions of the forced  Navier--Stokes equations in
  $\Omega$ (subject to either no-slip Dirichlet boundary condition,
   in the presence of physical boundary, or in the whole space
   or in a periodic box)
\begin{equation*}
\del_t u_\nu +\nabla\cdot (u_\nu \otimes u_\nu)-\nu\Delta u_\nu +
\nabla p_\nu= f\,.
\end{equation*}
Here $u_\nu$ and $p_\nu$ are random variables which depend on $(x,t)$. We
will drop, below, the explicit dependence on $\mu$ when this does
not cause any confusion.


Multiplying the Navier--Stokes equations  by $u_\nu(x,t,\mu)$ and
{\it assuming} that one has the following equality
\begin{eqnarray*} \frac 1
2\del_t |u_\nu(x,t,\mu)|^2 -\nabla_x \cdot(( \nu  \nabla_x u_\nu -
p_\nu
I) u_\nu)(x,t,\mu) + \\
\nu |\nabla_x u_\nu(x,t,\mu)|^2=  f(x,t) \cdot u_\nu(x,t,\mu)\,,
\end{eqnarray*}
we observe, that in the $2d$ case, the above relation is a proven
fact. However,  in the $3d$ case, the class of suitable solutions of
the Navier--Stokes equations, in the sense of
 Caffarelli--Kohn--Nirenberg,  are known to satisfy a weaker form of the
above relation, involving an inequality instead of equality (cf.
(\ref{suitsol})).

Thanks to the homogeneity assumption, the quantity
$$
\langle (( \nu  \nabla_x u_\nu - p_\nu I) u_\nu )(x,t,\cdot)\rangle
$$
does not depend on $x$, and therefore
$$
\langle \nabla_x \cdot ( (( \nu  \nabla_x u_\nu - p_\nu I) u_\nu
))(x,t,\cdot)\rangle= \nabla_x \cdot\langle (( \nu  \nabla_x u_\nu -
p_\nu I) u_\nu )(x,t,\cdot)\rangle =0\,.
$$
Thus, the averaged pointwise
energy relation
\begin{equation}
\frac 1 2\del_t\langle  |u_\nu(x,t,\cdot)|^2\rangle +\nu \langle
|\nabla u_\nu(x,t,\cdot)|^2\rangle = \langle f(x,t) \cdot
u_\nu(x,t,\cdot)\rangle \label{Karman}
\end{equation}
is obtained. The above  is often
called the Karman--Howarth relation, and  it implies that the
quantities
$$
\mathfrak{e}=\langle  |u_\nu(x,t,\cdot)|^2\rangle \,\,\hbox{
and }\,\, \epsilon(\nu)=\frac{1}{t}\int_0^t  \nu\langle |\nabla
u_\nu(x,s,\cdot)|^2\rangle ds
$$ are uniformly bounded in time, under reasonable assumptions
on the forcing term.

Finally, assume also that the random process is stationary in time.
Then, the equation (\ref{Karman}) gives an {\it a priori} estimate
for the   mean rate of dissipation of energy
$$
\epsilon(\nu)=\nu \langle  |\nabla u_\nu(x,t,\cdot)|^2\rangle\,,
$$
which is independent of $(x,t)$. Now, with the forcing term $f$
acting only on low Fourier modes, one may assume the existence  of a
region, called the inertial range, where the turbulence spectra
$E(|k|)$ behaves according to a universal power law. {\it The name
inertial refers to the fact that in this range of wave numbers, the
energy cascades from low modes to high modes with no leakage of
energy. That is, there is no viscous effects in this range and only
the inertial  term $(u\cdot \nabla)u$
 is active}. Combining the above hypotheses: homogeneity, isotropy
and stationarity of the random process, together with the  existence
of an inertial range of size
 $$
 A\le |k|\le B\nu^{-\frac12}\,,
 $$
 where the spectra behaves according to a power law,
one obtains finally, by a dimensional analysis and in the
three-dimensional case, the ``Kolmogorov law":
\begin{equation}
E(|k|)\simeq \epsilon(\nu)^{\frac 23 }|k|^{-\frac 53}  \label{kol1}\,.
\end{equation}
It is important to keep in mind the fact that this derivation is
based on the analysis of a random family of solutions. Therefore,
the formula (\ref{kol1}) combined with the formula (\ref{wigal})
implies that {\it in the average  }
the solutions have a spectra which behaves in
the turbulent regime according to the prescription.
\begin{eqnarray*}
&&\langle \hat u_\nu  (k,t,\cdot)\otimes \overline {\hat u_\nu
(k,t,\cdot)} \rangle = \frac{1}{(2\pi L)^d}  \int_{(\R/L \ZZ)^d}
e^{-i\frac {k \cdot r}L }
\langle u_\nu(y+\frac r 2,t,\cdot) \otimes u_\nu(y-\frac r2,t,\cdot)\rangle  dr\\
&&\simeq  \frac{\epsilon(\nu)^{\frac 23 }|k|^{-\frac 53}}{4\pi
 |k|^{d-1}} (I-\frac{k\otimes k}{|k|^2})\,.
\end{eqnarray*}

\begin{remark}

The main difficulty in the full justification of the above
derivation is the construction of  a probability measure, $dm$, on
the ensemble of solutions of the Navier--Stokes equations that
would satisfy the hypotheses of homogeneity, isotropy and
stationarity. In particular, the construction of such measure should
be uniformly valid when the viscosity $\nu$ tends to $0\,.$ See, for
instance, the books of Vishik and Fursikov \cite{Vishik-Fursikov}
and Foias et al. \cite{FMRT} for further study and references
regarding this challenging problem.

The next difficulty (which is a controversial subject) is the
justification for the spectra of an inertial range with a power law.
In Foias et al. \cite{FMRT} (see also Foias \cite{Foias})  it was
established, for example, the existence of an inertial range of wave
numbers where one has a forward energy cascade. However, we are
unaware of a rigorous justification  for a power law in this
inertial range.

Nevertheless, the construction of a power law of the spectra is
often used as a benchmark for validation of numerical computations
and experiments. Since, in general, one would have only one run of an
experiment, or one run of simulation, a Birkhoff theorem, which
corresponds to assuming an ergodicity hypothesis, is then used. This
allows for the replacement of the ensemble average by a time
average. For instance, one may assume, in presence of forcing term,
in addition to stationarity that {\it for almost all solutions},
i.e. almost every $\mu$, one has
\begin{eqnarray*}
\lim_{T\rightarrow \infty} {\frac 1 T} \int_0^T u_\nu(y+\frac r
2,t,\mu ) \otimes u_\nu(y-\frac r2,t,\mu )dt=\\
\langle u_\nu(y+\frac
r 2,t,\cdot) \otimes u_\nu(y-\frac r2,t,\cdot)\rangle \,,
 \end{eqnarray*}
 which would give the following relation, for almost every
 solution $u_\nu$,
 \begin{eqnarray*}
&&\lim_{T\rightarrow \infty} \frac 1 T\int_0^T u_\nu
(k,t,\mu)\otimes \overline {\hat u_\nu (k,t,\mu)}  dt = \\
&&\frac{1}{(2\pi )^d}  \int_{\R^d} e^{-i k\cdot r}
\langle u_\nu(y+\frac r 2,t,\cdot) \otimes u_\nu(y-\frac r2,t,\cdot)\rangle  dr \\
&&= \langle \hat u (k,t,\cdot)\otimes \overline {\hat u (k,t,\cdot)}
\rangle \simeq  \frac{\epsilon(\nu)^{\frac 23 }|k|^{-\frac 53}}{4\pi
 |k|^{d-1}} (I-\frac{k\otimes k}{|k|^2})\rangle\,.
\end{eqnarray*}
\end{remark}

\subsection{ Comparison between deterministic and statistical spectra.}

The deterministic point of view considers families of  solutions
$u_\nu$ of the Navier--Stokes dynamics, with viscosity $\nu>0$, and
interprets the notion of turbulence in terms of the weak limit
behavior (the asymptotic  behavior of such sequence as $\nu\to 0$)
with the Wigner spectra:
\begin{eqnarray*}
&& \hskip 0.5in \widehat{ RT } (x,t,dk)=\\
&&\lim_{\nu \rightarrow 0} \frac1 {(2\pi) ^d}\int_{\R^d_y}
\!\!e^{ik\cdot y} ((u_\nu - \overline u) (x-\frac {\sqrt {\nu}}2
y,t)\otimes (u_\nu-\overline u)(x+\frac {\sqrt {\nu}}2 y,t))dy.
 \end{eqnarray*}
As we have already observed earlier, this is a local object (it
takes into account the $(x,t)$ dependence). Moreover, one could
define the support of turbulence, for such family of solutions, as
the support of the measure ${\mathrm {Trace}}\widehat{ RT }
(x,t,dk)$. Of course, determining such a support is extremely
hard and is a configuration dependent problem. This is perfectly
described in the sentences of Leonardo da Vinci, who is   very
often quoted, and in particular in  page 112 of \cite{FR}:

\vskip 0.125in
 {\it \hskip 0.5in doue la turbolenza dellacqua sigenera

\hskip 0.5in doue la turbolenza dellacqua simantiene plugho

\hskip 0.5in  doue la turbolenza dellacqua siposa.}

\vskip 0.125in
 Up to this point nothing much can be said  without
 extra hypotheses, except that the formula (\ref{freq2})
indicates
the existence of an essential, if not an ``inertial", range
\begin{equation}
A\le |k|\le \frac B{\sqrt{\nu}}\,.\label{detinr}
\end{equation}
On the other hand, the statistical theory of turbulence starts from
the hypotheses, that seem  difficult  to formulate in a rigorous
mathematical setting, concerning the  existence of statistics (a probability
measure) with respect to which the  two points correlations
 for any family of solutions,  $u_\nu$, of the Navier--Stokes
 equations are  homogeneous and isotropic. Under these
 hypotheses, one proves  properties on the decay
 of the spectra of turbulence. Moreover,
with all these assumptions one obtains, by simple dimension
analysis, {\it for averages of solutions with respect to the
 probability measure $dm(\mu)$ }, the following formula
 in the inertial range:
\begin{equation}
\langle |\hat u_\nu(k,t,\cdot)|^2\rangle
=\frac{1}{(2\pi)^4}\frac{E(|k|)}{|k|^2}\simeq (\epsilon(\nu))^{\frac
23 }|k|^{-\frac {11}3}\label{kol2bis}\,.
 \end{equation}

 Finally, by assuming and using the stationarity (in time) of these
 two points correlations, and by applying, sometimes,  the
 Birkhoff ergodic theorem, one should be able to obtain,
 for almost every solution, the formula
 \begin{eqnarray*}
&&\lim_{T\rightarrow \infty} \frac 1 T\int_0^T (\hat u_\nu (k,t,\mu)\otimes \overline {\hat u_\nu (k,t,\mu)} ) dt =\nonumber\\
&& \langle \hat u_\nu (k,t,\cdot)\otimes \overline {\hat u_\nu
(k,t,\cdot)} \rangle \simeq  \frac{(\epsilon(\nu))^{\frac 23
}|k|^{-\frac 53}}{4\pi  |k|^{d-1}} (I-\frac{k\otimes
k}{|k|^2})\rangle\,.
\end{eqnarray*}
\begin{remark} Some further connections between these two aspects of spectra may be considered.

i) Assuming that near a point $(x_0,t_0)$ the Wigner spectra is
isotropic. Define the local mean dissipation rate of energy as
\begin{equation}
\epsilon_{(x_0,t_0)}(\nu)= \nu \int_0^\infty\!\int_\Omega
|\phi(x,t)|^2|\nabla u_\nu(x,t)|^2dxdt\,,
\end{equation}
with $\phi$ being a localized
function about $(x_0,t_0)$. Then one can prove that, for $|k|$ in
the range given by (\ref{detinr}),
\begin{eqnarray*}
 \frac1 {(2\pi)^d}\int_{\R^d_y} e^{ik\cdot y} ((u_\nu-\overline u) (x-\frac {\sqrt {\nu}}2 y,t)\otimes (u_\nu-\overline u)(x+\frac {\sqrt {\nu}}2 y,t))dy   \\
\sim \epsilon_{(x_0,t_0)}(\nu) |k|^{-\frac{11}3}(I-\frac{k\otimes
k}{|k|^2})\,,
 \end{eqnarray*}
  as   $\nu $ tends  to zero.

ii) Give sufficient conditions that will make the Wigner spectra
isotropic. This is in agreement with the fact that this spectra
involves only a small scale phenomena; thus, it is a reasonable
hypothesis. However, the example constructed by Cheverry \cite{CHE}
shows that this has no chance  of always being true.

iii) Another approach for establishing the existence of an inertial
range for forward energy cascade in $3d$, and forward enstrophy
cascade in $2d$, is presented in \cite{FMRT} (see also references
therein). This approach is based on the statistical stationary
solutions of the Navier--Stokes equations. These are time
independent probability measures which are invariant under the
solution operator of the Navier--Stokes equations. Furthermore, in
Foias \cite{Foias} some semi-rigorous arguments  are presented to
justify the Kolmogorov power law of the energy spectrum.

\end{remark}

\section{ Prandtl and Kelvin--Helmholtz problem}

In this section it is assumed that the sequence of solutions
$\{u_\nu\}$ of the Navier--Stokes equations with no-slip Dirichlet
boundary condition (in the presence of physical boundary) converges
to the solution of the Euler equations. According to   Theorem
\ref{kato} of Kato, in this situation one has
\begin{eqnarray}
\lim_{\nu\rightarrow 0} \nu \int_0^T\int_{\{x\in
\Omega:\,d(x,\del\Omega)<\nu\}}|\nabla
u_\nu(x,t)|^2dxdt=0\label{kb}\,.
\end{eqnarray}

However, since the tangential velocity of the solution of the Euler
equations is not zero on the boundary,  as $\nu\rightarrow 0$, a
boundary layer is going to appear. On the one hand, the scaling of
the boundary layer has to be compatible with the hypothesis
(\ref{kb}); and on the other hand the equations that model the
behavior in this boundary layer have to reflect the fact that the
problem is very unstable. This is because the instabilities (and
possible singularities) that occur near the boundary may not remain
confined near the boundary, and will in fact propagate inside the
domain by the nonlinear advection term of the Navier--Stokes
equations. These considerations explain why the Prandtl equations
(PE) of the boundary layer are complicated.

There are good reasons to compare the Prandtl equations with the
Kelvin--Helmholtz problem (KH):

1. Even though some essential issues remain unsolved for KH, it is
much better understood from the mathematical point of view than
the PE problem.  However, the two problems
share similar properties such as instabilities and appearance of
singularities.

2. At the level of modeling, in particular for the problem
concerning the wake behind an air plane and the vortices generated
by the tip of the wings, it is not clear if turbulence should be
described by singularities in KH or PE (or both)!

For the sake of simplicity, these problems are considered in the
$2d$ case, and for the PE in the half space $x_2>0\, .$

\subsection{ The Prandtl Boundary Layer}
One starts  with the $2d$ Navier--Stokes equations in the half plane
$x_2>0$, with the no-slip boundary condition $u^\nu(x_1,0,t)\equiv0$
{\begin{eqnarray}
&&\del_t u_1^\nu - \nu\Delta u_1^\nu +
u_1^\nu\del_{x_1} u_1^\nu +u_2^\nu\del_{x_2} u_1^\nu +\del_{x_1}
p^\nu=0\,,
\\
&&\del_t u_2^\nu - \nu \Delta u_2^\nu + u_1^\nu\del_{x_1} u_2 ^\nu+
u_2^\nu\del_{x_2} u_2^\nu +\del_{x_2} p^\nu=0\,,
\\
&&\del_{x_1} u_1^\nu+\del_{x_2} u_2^\nu =0 \\
&& u_1^\nu(x_1,0,t)=u_2^\nu(x_1,0,t)=0  \hbox { on }  x_1\in \R\,,
\end{eqnarray}
and assumes that inside the domain (away from the boundary) the
vector field $u^\nu(x_1,x_2,t)$ converges to the solution $u_{\rm
int} (x_1,x_2,t)$ of the Euler equations with the same initial data.
The tangential component of this solution on the boundary, $x_2=0$,
and of the pressure are  denoted by
$$U(x_1,t)=u_1^{\rm int} (x_1,0,t)\,, \quad \tilde P(x_1,t)=p(x_1,0,t)\,. $$
Then one introduces the scale  $\epsilon=\sqrt \nu$. Taking into
account that the  normal  component of the velocity remains $0$ on
the boundary, one uses the following ansatz, which corresponds to a
boundary layer  in a parabolic PDE problem
\begin{equation}
 \Bigg(\begin{array}{c}
 u_1^\nu(x_1,x_2) \\
u_2^\nu(x_1,x_2)
\end{array}  \Bigg)=\Bigg(\begin{array}
{c}
\tilde{u}_1^\nu(x_1,\frac{x_2}{\epsilon}) \\
\epsilon \tilde {u}_2^\nu(x_1,\frac{x_2}{\epsilon})
\end{array}\Bigg)+ u^{\mathrm {int}}(x_1,x_2)\label{couchlim}\,.
\end{equation}
Inserting the right hand side of (\ref{couchlim})  into the
Navier--Stokes equations, and returning to the notation $(x_1,x_2)$ for
the variables
$$
X_1=x_1, X_2=\frac{x_2}{\eps}\,,
$$
and letting $\epsilon$ go to zero, one obtains formally the equations:
\begin{eqnarray}
&&  \uT_1(x_1,0,t )+U_1(x_1,0,t)=0\,,  \label{Eulerext}\\
&& \del_{x_2} \tilde p(x_1,x_2)=0\Rightarrow  \tilde p  (x_1,x_2,t) =\tilde P (x_1,t)\label{Eulerext2}\\
&&\del_t \uT_1  - \del_{x_2}^2 \uT_1 + \uT_1 \del_{x_1} \uT_1
+\uT_2 \del_{x_2} \uT_1   =\del_{x_1} \tilde P(x_1,t)\,, \label{Eulerext3}\\
&&
\del_{x_1} \uT_1 +\del_{x_2} \uT_2  =0\,,\,\, \uT_1(x_1,0)=\uT_2(x_1,0) =0 \hbox { for } x_1\in \R\,,
\\
&&\lim _{x_2\rightarrow \infty} \uT_1(x_1,x_2)=\lim _{x_2\rightarrow
\infty} \uT_2(x_1,x_2)=0\,.
\end{eqnarray}
\begin{remark}
 As an indication of the validity of the Prandtl equations we
 observe that (\ref{couchlim}) is consistent with Kato
 Theorem \ref{kato}. Specifically, thanks to
 (\ref{couchlim}) one has
\begin{equation*}
\nu\int_0^T\!\!\!\int_{\Omega\cap \{d(x,\del \Omega)\le c \nu\}}
|\nabla u^\nu(x,s)|^2dxds\le C\sqrt{ \nu}\,.
\end{equation*}
\end{remark}
\begin{remark}
The following example, constructed by Grenier \cite{GR}, shows that
the Prandtl expansion cannot always be valid. In the case when the
solutions are considered in the domain
$$({\R_{x_1}}/{\mathbf Z})\times\R_{x_2}^+\,.$$
Grenier starts with a solutions $u^\nu_{ref}$ of the pressureless
Navier--Stokes equations given by
\begin{eqnarray*}
&&u^\nu_{ref}=(u_{ref}(t,y/\sqrt{\nu}),0)\\
&&\del_t u_{ref} -\del_{YY}u_{ref}=0\,,
\end{eqnarray*}
where $Y= y/\sqrt{\nu}$. Using a convenient and explicit choice of
the function $u_{ref}$, with some sharp results on instabilities, a
solution of the Euler equations of the form
\begin{equation*}
\tilde u =u_{ref} +\delta v + O(\delta^2 e^{2\lambda t}) \quad {\rm
for } \quad 0<t<\frac 1{\lambda \log \delta}
\end{equation*}
is constructed. It is then shown that the vorticity generated by the
boundary for the solution Navier--Stokes equations (with the same
initial data) is too strong to allow for the convergence of the
Prandtl expansion. One should observe, however,  that once again
this is an example which involves solutions with infinite energy. It
would be interesting to see if such an example could be modified to
belong to  the class of finite energy solutions; and then to analyze
how   the modified finite energy solution might violate the Kato
criteria mentioned in Theorem \ref{kato}.
\end{remark}
It is  important to observe that in their mathematical properties
the PE exhibit the pathology of the situation that they are trying
to model. First one can prove the following proposition.
\begin{prop} Let
$T>0$ be a finite positive time, and let $(U(x,t),P(x,t))\in
C^{2+\alpha}(\R_t\times (\R_{x_1}\times \R_{x_2}^+))$ be a smooth
solution of the $2d$ Euler equations satisfying, at time $t=0$, the
compatibility condition $U_1(x_1,0,t)=U_2(x_1,0,t)=0$ (notice that
only the boundary condition $U_2(x_1,0,t)=0$ is preserved by the
Euler dynamics). Then the following statements are equivalent:

i) With initial data $\uT(x,0)=0\,,$ the boundary condition
$\uT_1(x_1,0,t)=U_1(x_1,0,t)$ in (\ref{Eulerext}), the right-hand
side in (\ref{Eulerext2})  given by
$\tilde P(x_1,t)= P(x_1,0,t)\,,$
and the Prandtl equations  have a smooth solution $\uT(x,t)$, for
$0<t<T$\,.

ii) The solution $u_\nu(x,t)$, of the Navier--Stokes equations with
initial data $u_\nu(x,0)=U(x,0)$  and with no-slip boundary
condition at the boundary $x_1=0$, converges in $C^{2+\alpha}$ to
the solution of the Euler equations, as $\nu \to 0$.
\end{prop}
The fact that  statements  i) and ii) may be violated for some
$t$ is related to the appearance of a detachment zone, and the
generation of turbulence. This is well illustrated in the analysis
of the Prandtl equations   written in the following
simplified form
\begin{eqnarray}
&&\del_t \uT_1  - \del_{x_2}^2 \uT_1 + \uT_1 \del_{x_1} \uT_1 +\uT_2
\del_{x_2} \uT_1   =0\,, \label{pes1}
\del_{x_1} \uT_1 +\del_{x_2}
\uT_2  =0
\\
&&\uT_1(x_1,0)=\uT_2(x_1,0) =0 \hbox { for } x_1\in \R \label{pes2}
\\
&&\lim _{x_2\rightarrow \infty} \uT_1(x_1,x_2)=0, \\
&& \tilde u_1(x_1,x_2,0)=\tilde u_0(x_1,x_2)\,.
\end{eqnarray}
Regularity in the absence of detachment corresponds to a theorem of
Oleinik \cite{OL}. She proved that global smooth solutions of the
above system do exist provided the initial profile is monotonic,
i.e. for any initial profile satisfying
\begin{equation*}
\uT(x,0)=(\uT_1(x_1,x_2),0)\,,\quad \del_{x_2} \uT_1(x_1,x_2)
\not=0\,.
\end{equation*}
On  the other hand, initial conditions with ``recirculation
properties" leading to a finite time blow up have been constructed
by E and Engquist \cite{WEE} and \cite{WEE2}.  An interesting aspect of these
examples is that the blow up generally does not occur   on
the boundary, but rather inside the domain.

The above pathology appears in the fact that the PE
is highly unstable. This comes from the determination
$\uT_2$  in term of $\uT_1$ by the equation
\begin{equation*}
\del_{x_1} \uT_1 +\del_{x_2} \uT_2  =0\,.
\end{equation*}
Therefore, it is only with analytic initial data (in fact analytic
with respect to the tangential variable is enough) that one can
obtain (using an   abstract version of the Cauchy--Kowalewskya
theorem) the existence of a smooth solution of the Prandtl equation
for a finite time and the convergence to the solution of the Euler
equations during this same time  (Asano \cite{AS}, Caflisch-Sammartino \cite{CS},
and Cannone-Lombardo-Sammartino \cite{CLS}.)
\subsection{The Kelvin--Helmholtz problem}\label{ke}
The Kelvin--Helmholtz (KH) problem concerns the evolution of a
solution of the $2d$ Euler equations
\begin{equation}
\del_t u +\nabla \cdot (u\otimes u)+\nabla p =0 \,, \quad \nabla
\cdot u =0
\end{equation}
with initial vorticity $\omega(x,0)$ being a measure concentrated on
a curve $\Gamma(0)\,.$

This is already simpler than the PE because the pathology, if any,
should  in principle be concentrated on a curve. Furthermore, the
dynamics in this case inherits the general properties of the $2d$
dynamics. In particular,  it will obey the equation
\begin{equation*}
\del_t(\nabla \wedge u) +(u \cdot \nabla) (\nabla \wedge u)=0\,,
\end{equation*}
for the conservation (for smooth solutions) of the density of the
vorticity. Therefore, one  can guarantee the existence of a weak
solution when the initial vorticity, $\omega(x,0)$\,,  is a Radon
measure. This was done first by Delort, assuming that the initial
measure has a distinguished sign \cite{DE}. Then the result was
generalized to situations where the change of sign was simple enough
\cite{LNZ}. However, this remarkable positive result is impaired by
the non-uniqueness result of Shnirelman \cite{Sh}.

For smooth solutions of the KH, i.e. the ones with  vorticity $\omega$ - a
bounded Radon measure with support contained in curve
$\Sigma_t=\{r(\lambda,t),\, \lambda\in \R\}\,$ - the velocity field
is given, for $x\notin \Sigma_t$, by the so called Biot--Savart law
\begin{eqnarray}
&&u(x,t)=\frac1{2\pi}R_{\frac{\pi}{2}}  \int \frac {x-r'}
{|x-r'|^2} \omega(r',t)ds'\nonumber\\
&&:=\frac1{2\pi}R_{\frac{\pi}{2}} \int \frac {x-r(\lambda',t)}
{|x-r(\lambda',t)|^2} \omega(r(\lambda',t),t)\frac{\del
s(\lambda',t)}{\del\lambda'}d\lambda'\label{bios}\,.
\end{eqnarray}
Here  $r(\lambda,t)$, for $\lambda \in \R$, is a parametrization of
the curve $\Sigma_t$,  and $s(\lambda,t)= |r(\lambda,t)|$ is the
corresponding arc length.  $R_{\frac{\pi}{2}}$  denotes the
$\frac{\pi}{2}-$counterclockwise rotation matrix
$$
R_{\frac{\pi}{2}}=\left(
\begin{array}{cc}0& -1\\
1&0
\end{array} \right)\,.
$$
Furthermore, as  $x$ approaches the curve $\Sigma_t\,$ the velocity
field $u$ admits the two-sided limits $u_\pm\,.$ By virtue of the
incompressibility condition one has the continuity condition for the
normal component of the velocity field, i.e.
\begin{equation*}
u_-\cdot \vec n=u_+\cdot \vec n\,,
\end{equation*}
where hereafter $\vec \tau$ and $\vec n$ will denote the unit
tangent and unit normal vectors to the curve $\Sigma_t$,
respectively. In addition, the average
$$\langle u \rangle=\frac{ u_++u_-}2$$
is given by the principal value of the  singular
integral appearing in (\ref{bios})
\begin{equation}
v=\langle u \rangle=\frac1{2\pi}R_{\frac{\pi}{2}}\, {\rm p.v.}
\!\!\int \frac {x-r'} {|x-r'|^2} \omega(r',t)ds'\,.
\label{average-velocity}
\end{equation}
Using the calculus of distributions one can show that, as long as
the curve $\Sigma_t$ is smooth,    the  velocity
field $u$, defined above, being a weak solution of the Euler
equations
\begin{equation*}
 \del_t u  + \nabla \cdot ( u\otimes u) + \nabla p=0\,,\,\,
 \nabla \cdot u=0\,,
\end{equation*}
 is equivalent to the  vorticity  density $\omega$
 and the curve $\Sigma_t$
 satisfying the coupled system of equations
 \begin{eqnarray}
&&\omega_t-\del_s\Bigg(\omega(\del_t r-v)\cdot \vec \tau\Bigg)=0\,,\label{rs}\\
&&(r_t-v)\cdot \vec n =0\,,\label{para}\\
&&v(r,t)=\frac{1}{2\pi} R_{\frac{\pi}{2}}\, {\rm p.v.} \!\!\int
\frac{r-r'}{|r-r'|^2}\omega(r',t)ds'\,.\label{bios2}
 \end{eqnarray}
The equations (\ref{rs}), (\ref{para}) and  (\ref{bios2}) do  not
completely determine   $r(\lambda,t)\,.$ This  is due to the freedom
in the choice of the parametrization of the curve $\Sigma_t$.
Assuming that $\omega\not=0$ one introduces a new parametrization
$\lambda(t,s)$ which reduces the problem to the equation
\begin{equation}
\del_t r (\lambda,t)  =\frac{1}{2\pi} R_{\frac{\pi}{2}}\, {\rm p.v.}
\!\!\int \frac{r(\lambda,t)
-r(\lambda',t)}{|r(\lambda,t)-r(\lambda',t)|^2}
d\lambda'\label{birot2}\,,
\end{equation}
or with the introduction of the complex variable  $z=r_1 +i r_2$,
where $r=(r_1, r_2)$, one obtains the Birkhoff--Rott equation
\begin{equation}
\del_t \overline z (\lambda,t) = \frac{1}{2\pi i } \,{\rm p.v.} \!\!
\int \frac{d\lambda'}{z(\lambda,t)-z(\lambda',t) }\label{birot3}\,.
\end{equation}

\begin{remark}
The following are certain mathematical similarities of the KH
problem with the PE:

1. As  for the PE one has for the evolution equation (\ref{birot3})
a local, in time, existence and uniqueness result in the class of
analytic initial data. This is done by implementing a version
of the Cauchy--Kowalewskya theorem (Bardos, Frisch, Sulem and Sulem
\cite{BFSS}.)

2.  As for the PE one can construct  solutions that blow up in
finite time.

3. One observes  that the singular  behavior in the experiments and
numerical simulations with the KH problem is  very similar   to the
one that is generated by the no-slip boundary condition when the
viscosity is approaching  zero.
\end{remark}

The best way to understand the structure of the KH is to use the
fact that the Euler  equations are invariant under both space and
 time translations, and under space rotations,  and to
 consider a weak solution of the $2d$ Euler
 dynamics either in the whole plane $\R^2$, satisfying
 $$
 u\in C((-T,T); L^2(\R^2))\,,\,\, T>0
 $$
 or subject to periodic boundary conditions satisfying
$$
u\in C((-T,T); L^2((\R/{\mathbf Z})^2)).
$$
 Assuming
 that in a small neighborhood $\U$ of the point
 $(t=0,z=0)$ the
vorticity is concentrated on a smooth curve in the complex plane
which takes the form:
\begin{equation}
z(\lambda,t)\!=  \!(\alpha t+\beta(\lambda\! +\!\eps f(\lambda,t))
\,,\,  f(0,0)\!=\!\nabla f(0,0)\!=\!0\,. \label{local2}
\end{equation}
Then using the relations $\nabla\cdot u =0\,,\, \nabla \wedge
u=\omega$ and the Biot--Savart law, one obtains
\begin{eqnarray}
&&\hskip -0.25in  \epsilon|\beta|^2  \del_t \overline f (\lambda,t)  = {\nonumber} \\
&&\frac{1}{2\pi i}\,{\rm p.v.} \!\!\int_{\{z(t,\lambda')\in \UU\}}
\frac{ d\lambda'}{ (\lambda-\lambda')
(1-\epsilon\frac{f(\lambda,t)-f(\lambda',t)}{\lambda-\lambda'})}\!+\!
E(z(\lambda,t))\label{parpro1}
\end{eqnarray}
where here, and in the sequel,  $E(z)$ denotes the ``remainder",
which is  analytic with respect to $z$. Next,  use the expansion
\begin{eqnarray}
&&\frac1{2\pi}\, {\rm p.v.} \!\!\int \frac
{d\lambda'}{(\lambda-\lambda')(1-\epsilon\frac{f(\lambda,t)-f(\lambda',t)}{\lambda-\lambda'})}d\lambda'
= \label{parpro2}\\
&& \frac{\epsilon}{2\pi}\, {\rm p.v.}
\!\!\int\frac{f(\lambda,t)-f(\lambda',t)}{(\lambda-\lambda')^2}d\lambda'
+\sum_{n\ge 2}\frac{\eps^n}{2\pi}\, {\rm p.v.} \!\!
\int\frac{(f(\lambda,t)-f(\lambda',t))^n}{(\lambda-\lambda')^{(n+1)}}d\lambda'\,,\nn
\end{eqnarray}
and implement the following  formulas concerning the Hilbert
transform
\begin{eqnarray}
&&\frac1{2\pi} \, {\rm p.v.} \!\!\int\frac{f(\lambda,t)-f(\lambda',t)}{\lambda-\lambda'}d\lambda'= -\frac{i}{2} \, {\mathrm{sgn}}(D_\lambda)f\,,\\
&&\,\, \frac1{2\pi} \, {\rm p.v.}
\!\!\int\frac{f(\lambda,t)-f(\lambda',t)}{(\lambda-\lambda')^2}d\lambda'=|D_\lambda|f
\end{eqnarray}
to deduce, from (\ref{parpro1}) and (\ref{parpro2}), that the real
and imaginary parts of $f(\lambda,t)=X(\lambda,t)+iY(\lambda,t)$ are
local solutions of the system
{\begin{eqnarray}
&&\del_tX  =\frac{1}{2|\beta|^2}|D_\lambda| Y+\epsilon R_1(X,Y)+E_1(X,Y) \\
&&\del_tY =\frac{1}{2|\beta|^2}|D_\lambda|X+ \epsilon R_2(X,Y) +
E_2(X,Y)
\end{eqnarray}}
or
{\begin{eqnarray}
&&(\del_t^2 +\frac{1}{4|\beta|^4} \del_\lambda^2)Y=\epsilon (\del_tR_1(X,Y)-\frac{1}{2|\beta|^2}|D_\lambda|R_2(X,Y))\nonumber\\
&&+\del_t E_1(X,Y)-\frac{1}{2|\beta|^2}|D_\lambda| E_2(X,Y) \label{el1}\\
&&(\del_t^2  +\frac{1}{4|\beta|^4} \del_\lambda^2)Y =\epsilon (\del_tR_2(X,Y)-\frac{1}{2|\beta|^2}|D_\lambda|R_1(X,Y))\nonumber\\
&&+\del_t E_2(X,Y)-\frac{1}{2|\beta|^2}|D_\lambda|
E_1(X,Y)\, .\label{el2}
\end{eqnarray}}
In (\ref{el1} ) and (\ref{el2}) the terms
\[\del_t E_1(X,Y)-\frac{1}{2|\beta|^2}|D_\lambda| E_2(X,Y) \hbox{ and }
\del_t E_2(X,Y)-\frac{1}{2|\beta|^2}|D_\lambda| E_1(X,Y)
\]
are the first order derivatives of analytic functions with respect
to $(X,Y)$\, while the terms
\[\del_tR_1(X,Y)-\frac{1}{2|\beta|^2}|D_\lambda|R_2(X,Y)) \hbox{ and }
(\del_tR_2(X,Y)-\frac{1}{2|\beta|^2}|D_\lambda|R_1(X,Y))
\]
are the second order derivative of analytic functions with a small
$\epsilon$ prefactor. Therefore, one observes that, up to a
perturbation, the KH problem behaves like a second order constant
coefficient elliptic equation. This fact has several important
consequence.

1. It explains why the evolution equation is well-posed only for a
short time, and with initial data that belongs to the class of
analytic functions. It is like solving an elliptic equation
simultaneously with both the  Dirichlet and Neumann boundary
conditions.

2. It is a tool for the construction of the solutions that blow up
in finite time.

3. It explains, by an indirect regularity argument,
 the very singular
behavior of the solution after the first break down of its
regularity.

These three points are discussed in further details below.

\subsubsection{Local solution}
When the curve $\Sigma_t$ is a graph of a function, say $y=y(x,t)$,
 the equations (\ref{rs}) and (\ref{para}) become
\begin{eqnarray}
&&y_t+y_x v_1=v_2\,,\quad
\del_t \omega+\del_x(v_1 \omega ) =0 \label{xpar1}\\
&&v_1(t,x)=-\frac{1}{2\pi}\, {\rm p.v.} \!\!\int_\R\frac{y(x,t)-y(x',t)\omega(x',t)}{(x-x')^2+(y(x,t)-y(x',t))^2}dx' \label{xpar2}\\
&&v_2(t,x)=\frac{1}{2\pi}\, {\rm p.v.}
\!\!\int_\R\frac{(x-x')\omega(x',t)}{(x-x')^2+(y(x,t)-y(x',t))^2}dx'\label{xpar3}\,,\end{eqnarray}
where   $(v_1,v_2)= v$ is the average velocity given in
(\ref{average-velocity}). Therefore, the above evolution equations
involve two unknowns $y(x,t)$ and $\omega (x,t)$, which is also the
case for the Birkhoff--Rott equation (\ref{birot3}), where the two
unknowns are the two components of $r(s,t)=(x(s,t),y(s,t))$, or of
$z(s,t)=x(s,t)+iy(s,t)$. In fact, since the Birkhoff--Rott equation
has been obtained by choosing the density of vorticity as a parameter,
one recovers this vorticity by the formula
$$
\omega(s,t)= \frac1{|\del_s z(s,t)|}\,.
$$
Since the system is a local perturbation of a second order  elliptic
equation then imposing two constraints at $t=0$ is similar to
solving this elliptic equation simultaneously with both Neuman and
Dirichlet boundary conditions. It is known that in the absence of
stringent compatibility conditions (they are related by the so called
Dirichlet to Neumann operator) such a problem can be solved only
locally and with analytic data. This is the reason   why the
solution of (\ref{xpar1})--(\ref {xpar3}) is obtained  locally in
time under the assumption that the functions $y(x,0) $ and
$\omega(x,0)$ are analytic.
\subsubsection{Singularities} For the construction of
singularities one follows the same  idea, and furthermore, uses the
time reversibility of the $2d$ Euler equations. More precisely, if
one constructs solutions which are singular at $t=0$ and are regular
on the interval $(0,T]$ this will imply, just by changing the time
variable $t$ into $T-t$, the existence of smooth solutions at $t=0$
that blow up at $t=T$. The first result in this direction was
obtained by Duchon and Robert \cite{DR}; the initial condition on
the vorticity at $t=0$ is relaxed, and one assumes that the solution
$y(x,t)$ goes to zero as $t\rightarrow \infty$. Then one can
consider the system (\ref{xpar1})--(\ref {xpar3}) as a two point Dirichlet
boundary-value problem with $y(x,t)$  is given for $t=0$ and is
required to tend to zero as $t\rightarrow \infty$. Then by a perturbation method one
proves the following proposition.
\begin{prop}
There exists   $\epsilon>0$ such that for any initial data that
satisfies the estimate
\begin{equation}
y(x,0)=\int e^{ix\xi}  g(\xi)d\xi\,, \,
{\text{with}}\,\int|g(\xi)|d\xi\le \epsilon\,, \label{varb}
\end{equation}
the problem (\ref{xpar1})$-$(\ref {xpar3}) has a unique solution
which goes to zero as $t\rightarrow \infty\,. $ Furthermore, this
solution is analytic (with respect  to $(x,t)$) for all $t>0\,.$
\end{prop}
As  mentioned above,  this is a result about singularity
formations. It exhibits (by changing the time variable $t$ into
$T-t$) an example of solutions which are analytic at some time, but
with no more regularity at a later time than  what is allowed by the
equation (\ref{varb}). In fact, it was observed in some numerical
experiments \cite{M} and \cite{MBO} that the first break down of
regularity appears as a cusp on the curve $r(\lambda,t)\,.$ This
motivated  Caflsich and Orellana \cite{CO} to introduce the
function
\begin{equation}
f_0(\lambda,t)=(1-i)\{(1-e^{-\frac t
2-i\lambda})^{1+\sigma}-(1-e^{-\frac t 2+i\lambda})^{1+\sigma}\}\,,
\end{equation}
which  enjoys the following properties

i) For any $t>0$ the mapping $\lambda\mapsto f(\lambda,t)$ is
analytic.

ii) For $t=0$, the mapping $\lambda\mapsto f(0,\lambda)$ does not
belong to the H\"older space $C^{1+\sigma}$,  but it belongs to
every H\"older space $C^{1+\sigma'}$ with $0<\sigma'<\sigma\,.$

iii) The function
$$
z_0(\lambda,t)= \lambda +\epsilon f_0(\lambda,t)
$$
is an exact solution of the linearized Birkhoff--Rott
equation. More precisely, one has
\begin{equation}
\del_t \overline{f(\lambda,t)}= \frac{1}{2\pi}\, {\rm p.v.}
\!\!\int\frac{f(\lambda,t)-f(\lambda',t)}{(\lambda-\lambda')^2}d\lambda'
\,.
\end{equation}
Therefore, by using the ellipticity of this linear operator, one
can prove  by a
perturbation method the following proposition.
\begin{prop} \label{coo}For $\epsilon>0$, small enough,
there exist a function $r_\epsilon(\lambda,t) $ with the following
properties:

i) The function $\lambda\mapsto r_\epsilon(\lambda,t)$ is analytic
for $t>0\,.$

ii) The function $\lambda
+\epsilon(f_0(\lambda,t)+r_\epsilon(\lambda,t))$ is a solution of
the Birkhoff-Rott equation (\ref{birot3})\,.

ii) The function $\lambda \mapsto r_\epsilon(\lambda,t)$ is  (for $
\lambda\in \R \,,\, t\in \R_+$) uniformly bounded in $C^2\,. $
\end{prop}
As a consequence  of  Proposition  \ref{coo} (and of the
reversibility in time) one can establish the existence of analytic
solutions to the Birkhoff--Rott equation (\ref{birot3}), say in the
interval $0\le t<T$, such that at   time $t=T$ the map
$\lambda\mapsto z(\lambda,t)$ does not belongs to $C^{1+\sigma} $ at
the point $\lambda=0\,.$

\subsubsection{Analyticity and pathologic behavior after the break down of regularity}

The local reduction of the KH to the equation
\begin{eqnarray}
\epsilon|\beta|^2  \del_t \overline f (\lambda,t) = \frac{1}{2\pi
i}\,{\rm p.v.}\!\!\int_{z(\lambda',t)\in \UU}
\frac{d\lambda'}{ (\lambda-\lambda')
(1-\epsilon\frac{f(\lambda,t)-f(\lambda',t)}{\lambda-\lambda'})}\!+\!
E(z(\lambda,t))\label{parpro3}
\end{eqnarray}
requires obviously {\it some hypotheses} on the regularity of the
function $z(\lambda,t)$ near the point $(0,0)$. However, when this
reduction is valid it will, thanks to the ellipticity, imply that
the solution is $C^\infty$, and even analytic. Therefore, there
 appears to be a threshold (say {\bf T}) in the behavior of the solutions of the
KH. Existence of such {\it regularity threshold } is common in the
study of free boundary problems. This threshold is characterized by
the fact that any function having a regularity stronger than {\bf T}
is in fact analytic, and that there may exist solutions with less
regularity than {\bf T}. This has the following, practical, important
consequence: regularity of the solutions that are smooth for $t<T$
and singular after the time $t=T$ cannot be extended for $t\ge T$ by
solutions which are more regular than the threshold {\bf T}.
Otherwise, the above theorem would lead to a contradiction. This
fact explains why after the break down of regularity, the solution
becomes very singular.

For instance, it was shown by Lebeau \cite{LE} (and Kamotski and
Lebeau \cite{KL} for the local version) that any solution that is
near a point belongs to $C^\sigma_t(C^{1+\sigma}_\lambda)$ must be
analytic. As a consequence, if a solution constructed (by changing
the variable $t$ into $T-t$) according to the method of  Caflisch
and Orellana could be continued after time $t=T$, it would not be
in any H\"older space $C^{1+\sigma'}\,.$

Therefore, the  challenge (and an open problem) is the determination of
this  threshold of regularity that will imply analyticity. Up to
now, the best (to the best of our knowledge) known result is  due to
S.~Wu  ~\cite{Wu1} ~\cite{Wu2}. The hypothesis
$C_{loc}^\alpha(\R_t;C_{loc}^{1+\beta}(\R_\lambda))$ is replaced by
$H_{loc}^1( \R_t\times \R_\lambda)\,.$ The estimates are done by
explicitly using theorems of G.~David \cite{DA} saying that for all
chord arc curves $\Gamma: s\mapsto \xi(s),$  parameterized
by their arc length,  the Cauchy integral
operator
$$
C_\Gamma(f)= {\rm p.v.}\!\! \int
\frac{f(s')}{\xi(s)-\xi(s')}d\xi(s')
$$
is bounded in $L^2(ds)\,.$

The importance of this improvement is justified by the numerical
experiment of \cite{KR}. It is interesting to notice that these
results will apply to logarithmic spirals $r=e^\theta, \, \theta
\in\R$, but not to infinite length algebraic spirals. What is
observed is that from  cups singularity the solution evolves into a
spiral which behaves like an algebraic spiral, and therefore has an
infinite length. {\it The results of \cite{WU} provide an
explanation of the fact that the spiral has to be of an  infinite
length.}

After the appearance of the first singularity the solution becomes
very irregular. This leads to the issue of the definition of {\it
weak solutions} (solutions which are less regular than the threshold
{\bf T} ) {\it not of the Euler equations themselves, but of the
Birkhoff--Rott equation.} For instance, S.~Wu \cite{Wu1} and
\cite{Wu2} proposed the  following definition:

\vskip 0.125in

\noindent A weak solution is a function from $\R$ into $\mathbf{C},
\alpha\mapsto z(\alpha,t)$, for which the
following relation holds
$$\del_t(\int \overline z(\alpha,t)\eta(\alpha) d\alpha)=\frac{1}{4\pi i}
\int\!\!\!\int
\frac{\eta(\alpha)-\eta(\beta)}{z(\alpha,t)-z(\beta,t)}d\alpha
d\beta\,,
$$
for every  $\eta\in C^\infty_0(\R)$\,.

\vskip 0.125in

\noindent However, the problem is basically open because we have no
theorem concerning the existence of such a solution. Furthermore, for
physical reasons  weak solutions of the Birkhoff--Rott equation
should provide weak solutions of the incompressible Euler equations,
and in fact this is not always the case as it is  illustrated by
(cf. \cite{LNS}) the Prandtl--Munk example: Start from the vortex
sheet
 \begin{equation}
 \omega_0(x_1,x_2)=\frac{x_1}{\sqrt{1-x_1^2}}(\chi_{(-1,1)}(x_1)\otimes
 \delta(x_2))
 \end{equation}
 where $\chi_{(-1,1)}$ is  the characteristic function of
  the interval $(-1,1).$ By virtue of the the Biot--Savart
  law,  the velocity $v$ is constant
 \begin{equation}
  v  =(0,-\frac12)\,.
 \end{equation}
 The solution of the Birkhoff--Rott equation is given by the formula
 \begin{equation*}
 x_1(t)= x_1(0)\,, \,\,\, x_2(t)= \frac t2 \,,
 \,\,\, \omega(x_1,x_2,t  )=\omega_0(x_1,x_2+\frac t2)\,.
 \end{equation*}
On the other hand, it was observed in~\cite{LNS}  that the velocity
$u$ associated with this vorticity is {\it not} even a weak solution
of the Euler equations. In fact, one has
 \begin{equation}
 \nabla\cdot u=0\,\quad \hbox{ and } \del_t u + \nabla_x\cdot (u\otimes u)+\nabla p
 =F\,,
 \end{equation}
 where $F$ is given by the formula
 \begin{equation}
  F=\frac\pi 8  \left ((\delta(x_1+1,x_2+\frac t2)-\delta(x_1-1,x_2+\frac
  t2)),0\right )
 \,.\end{equation}
This has led  Lopes, Nussenweig and Sochet ~\cite{LNS} to propose a
weaker definition, which contains more freedom with respect to the
parameter, and may be more adapted.


\vskip0.2cm
{\it Acknowledgments}

The first author thanks the organizers of the Workshop on
Mathematical Hydrodynamics held in Moscow  June 2006 where he had
the opportunity to lecture on the subject. The authors are thankful
to the kind hospitalities of the  the Wolfgang Pauli Institute
(WPI), Vienna, and the  Bernoulli Center of the \'{E}cole
Polytechnique F\'{e}d\'{e}ral de Lausanne  where part of this work
was completed. This work was supported in part by the NSF grant
no.~DMS-0504619, the BSF grant no.~2004271,  and the ISF grant
no.~120/06.

\end{document}